\documentclass[11pt]{article}
\usepackage{epsfig}
\usepackage{amssymb}
\usepackage{amscd}
\usepackage{epsf}
\usepackage{amsfonts,amssymb,amsthm,amsmath}
 \usepackage{graphicx}
\usepackage[matrix,arrow,curve]{xy}
\usepackage{euscript}
\usepackage[sc]{titlesec}

\input cyracc.def

\long\def\comment#1\endcomment{\relax}

\newcounter{subsubsubsection}
\newcounter{subsubsubsubsection}

\makeatletter \@addtoreset{subsubsubsection}{subsubsection}
\@addtoreset{subsubsubsubsection}{subsubsubsection}

\makeatother

\newcommand{\sevafigc}[4]{\begin{figure}[h]\centerline{
 \epsfig{file=#1,width=#2,angle=#3}}
\bigskip\caption{#4}\end{figure}}

\comment
\endcomment

\newcommand{\udot}{{\:\raisebox{3pt}{\text{\circle*{1.5}}}}}
\DeclareMathSizes{11.1}{10}{8}{6}

\comment
\endcomment

\long\def\comment#1\endcomment{}

\theoremstyle{plain}

\newtheorem{theorem}{\sc Theorem}[section]
\newtheorem{lemma}[theorem]{\sc Lemma}
\newtheorem{conjecture}[theorem]{\sc Conjecture}
\newtheorem{prop}[theorem]{\sc Proposition}
\newtheorem{coroll}[theorem]{\sc Corollary}
\newtheorem{klemma}[theorem]{\sc Key-lemma}

\newtheorem{MainTheor}[theorem]{\sc Main Theorem}

\comment

\endcomment

\numberwithin{equation}{section}

\theoremstyle{plain}

\newtheorem{defn}[theorem]{\sc Definition}

\theoremstyle{exercise}
\newtheorem{remark}[theorem]{\sc Remark}
\newtheorem{remarks}[theorem]{\sc Remarks}
\newtheorem{example}[theorem]{\sc Example}

\newtheorem{notations}[theorem]{\sc Notations}

\DeclareMathOperator{\Hom}{Hom}

\comment

\newtheorem{theorem}{Theorem}

\newtheorem*{theorem*}{Theorem}

\newtheorem*{MainTheor}{Main Theorem}

\newtheorem*{lemma}{Lemma}

\newtheorem*{klemma}{Key-Lemma}

\theoremstyle{remark}

\newtheorem*{remark}{Remark}

\newtheorem*{example}{Example}

\theoremstyle{definition}

\endcomment

\newcommand{\eqto}{\mathrel{\stackrel{\sim}{\to}}}

\DeclareMathOperator{\Tr}{Tr}

\DeclareMathOperator{\ad}{ad}

\newcommand{\hdot}{{\:\protect\raisebox{3pt}{\text{\protect\circle*{1.5}}}}}

\newcommand{\mb}{\hdot}

\newcommand\End{\mathrm{End}}

\newcommand\Lie{\mathrm{Lie}}

\newcommand{\Ext}{\mathrm{Ext}}

\newcommand{\Conf}{\mathrm{Conf}}
\newcommand{\gl}{\mathfrak{gl}}

\newcommand{\Aff}{\mathrm{Aff}}

\newcommand{\U}{\mathcal{U}}

\newcommand{\g}{\mathfrak{g}}

\newcommand{\fin}{\mathrm{fin}}

\newcommand{\Comm}{\mathrm{Comm}}

\newcommand{\Der}{\mathrm{Der}}
\newcommand{\Inn}{\mathrm{Inn}}

\newcommand{\Hoch}{\mathrm{Hoch}}
\newcommand{\poly}{\mathrm{poly}}

\newcommand{\tot}{\mathrm{tot}}
\newcommand{\F}{\mathcal{F}}

\newcommand{\Assoc}{\mathrm{Assoc}}
\newcommand{\Cobar}{\mathrm{Cobar}}

\renewcommand{\k}{\Bbbk}
\newcommand{\Pois}{\mathrm{Pois}}
\newcommand{\aux}{\mathrm{aux}}
\newcommand{\Star}{\mathrm{Star}}
\newcommand{\loc}{\mathrm{loc}}
\newcommand{\Vect}{\mathrm{Vect}}
\newcommand{\CE}{\mathrm{CE}}

\title{ \sc An $L_\infty$ algebra structure on polyvector fields}
\author{\sc Boris Shoikhet}
\date{}

\begin{document}\maketitle

\begin{abstract}
\comment
In this paper we construct an $L_\infty$ structure on polyvector
fields on a vector space $V$ over $\mathbb{C}$ where $V$ may be
infinite-dimensional. We prove that the constructed $L_\infty$
algebra of polyvector fields is $L_\infty$ equivalent to the
Hochschild complex of polynomial functions on $V$, even in the
infinite-dimensional case. For a finite-dimensional space $V$, our
$L_\infty$ algebra is equivalent to the classical Schouten-Nijenhuis
Lie algebra of polyvector fields. For an infinite-dimensional $V$,
it is essentially different. In particular, we get the higher
obstructions for deformation quantization in infinite-dimensional
case.
\endcomment
It is well-known that the Kontsevich formality [K97] for Hochschild cochains of the polynomial algebra $A=S(V^*)$ fails if the vector space $V$ is infinite-dimensional. In the present paper, we study the corresponding obstructions. We construct an $L_\infty$ structure on polyvector fields on $V$ having the even degree Taylor components.  The degree 2 component is given by the Schouten-Nijenhuis bracket, but all its higher even degree components are non-zero. We prove that this $L_\infty$ algebra is quasi-isomorphic to the corresponding Hochschild cochain complex. We prove that our $L_\infty$ algebra is $L_\infty$ quasi-isomorphic to the Lie algebra of polyvector fields on $V$ with the Schouten-Nijenhuis bracket, if $V$ is finite-dimensional.
\end{abstract}

\section*{Introduction}
\subsection{}
The formality theorem of Maxim Kontsevich [K97] is one of the most important breakthroughs in the deformation theory of algebraic structures. It says that the dg Lie algebra of Hochschild cochains (with the Gerstenhaber Lie bracket) is quasi-isomorphic, as a dg Lie algebra, to its cohomology,
for an arbitrary regular commutative algebra of finite type over $\mathbb{C}$. The statement was firstly proven for the polynomial algebra $A=S(V^*)$, for a finite-dimensional vector space $V$ over $\mathbb{C}$, by making use a kind of Feynmann diagram expansion in topological quantum field theory. Later on, another proof was found by Dmitri Tamarkin  [T], using some operadic methods. For both proofs, the assumption that $V$ is {\it finite-dimensional}, is essential.\footnote{See Section \ref{sectionf} below for an explanation why the assumption $\dim V<\infty$ is essential for Kontsevich's proof, and Appendix B for the failure of Tamarkin's proof when $\dim V=\infty$.}  In Appendix A, we provide a proof of the general statement of the failure of the Kontsevich formality for an infinite-dimensional $V$.\footnote{In general, by the Kontsevich formality for $S(V^*)$ we always mean $\gl(V)$-equivariant formality. Note that our proof in Appendix A needs only a weaker equivariance.} 

This paper grew up from the author's attempt to single out the obstructions to the Kontsevich formality in the infinite-dimensional case.

Recall that the Kontsevich formality (for the case of polynomial algebra $S(V^*)$) is equivalent to existence of an $L_\infty$ quasi-isomorphism 
$$
\mathcal{L}\colon T_\poly(V)\to\Hoch^\udot (S(V^*))
$$
When $V$ is infinite-dimensional, both sides can be defined in an appropriate way (see Section \ref{section1.0}-\ref{section1.2}, they are denoted by $T_\fin(V)$ and $\Hoch_\fin^\udot(S(V^*))$, correspondingly). We show that, as complexes, the l.h.s. and the r.h.s. are quasi-isomorphic, in Theorem \ref{theoremhkr} (that is, ``the Hochschild-Kostant-Rosenberg theorem holds''). However, any $L_\infty$ quasi-isomorphisms {\it fails to exist}. 

In this paper, we construct a new $L_\infty$ structure on $T_\fin^\mathcal{L}(V)$ whose underlying graded vector space is $T_\fin(V)$, called ``the exotic $L_\infty$ structure''. Its Taylor components $\mathcal{L}_2,\mathcal{L}_4,\mathcal{L}_6,\dots$ are all nonzero, with $\mathcal{L}_2$ equal to the Schouten-Nijenhuis bracket (the odd degree components vanish by a symmetry reason). Our Main Theorem \ref{theoremmain} states that there is an $L_\infty$ quasi-isomorphism 
$$
\mathcal{L}_\fin\colon T_\fin^\mathcal{L}(V)\to \Hoch^\udot_\fin(S(V^*))
$$
regardless of is $V$ finite- or infinite-dimensional. Then the Kontsevich formality theorem implies that, for a finite-dimensional $V$, the $L_\infty$ algebras $T_\poly(V)$ (with the Schouten-Nijenhuis bracket and vanishing higher Taylor components) and $T_\fin^\mathcal{L}(V)$ are $L_\infty$ quasi-isomorphic.

\subsection{} 
The previous archive version(s) of this paper date(s) back to 2008. 
Recently, some other papers making use and further developing the ideas and results of this paper have appeared, see e.g. [KMW], [MW1], [MW2].
The author believes that the constructions introduced in this paper may find more fruitful applications in near future.

This 2017 version is essentially improved and expanded. The most important changes are:
(i) Section 1 is mostly rewritten with several new proofs and examples; (ii) we supplied the paper with two Appendices, A and B.  Appendix A provides a general proof of the failure on the Kontsevich formality for $S(V^*)$ for an infinite-dimensional $V$, and Appendix B shows why the Tamarkin proof [T] does not work in the infinite-dimensional case.
As well, we corrected English (which hopefully has improved since 2008).

\subsection*{}
\subsubsection*{Acknowledgements}
I am indebted to Boris Feigin who shared with me, around
'98-'99, some his conjecture on infinite-dimensional Duflo formula, which stimulated my work on 
infinite-dimensional formality in general. 

Maxim Kontsevich found a 
mistake in the proof of Lemma 2.2.4 of the first archive version, which failed Theorem
2.2.4 therein. The problem had been fixed fixed in a later archive version, by an interpretation of our former integrals as
the Taylor components of a new $L_\infty$ structure on polyvector fields, see Section 2. As well, the proof of Lemma \ref{lemmaa1} in Appendix A was communicated to the author by Maxim Kontsevich. I am
thankful to Maxim for his interest in my work and for the fruitful
correspondence. 

I am thankful to Pavel Etingof for his interest and
for his many suggestions.

The author is thankful to the anonymous referee for his careful reading of the paper, and for his remarks and suggestions which helped to improve the exposition.

The biggest part of the paper was completed during the author's 5-year appointment at the University of Luxembourg, when he was a member of the research group of Prof. Martin Schlichenmaier. I am thankful to Martin Schilechenmaier for his kindness and his personal participation, as well as for very nice working atmosphere, which made this paper possible to come into existence. 

The work was partially supported by the research grant R1F105L15 of the
University of Luxembourg and by the research grant nr. 6525 ``Kredieten aan Navorsers'' of Flemish Research Foundation (FWO).

\section{The set-up}\label{section1}
\subsection{}\label{section1.0}
It is known that the Kontsevich formality theorem [K97] for the Hochschild cochains of the polynomial algebra $A=S(V^*)$ fails, when the dimension of $V$ is infinite (we show in Section \ref{sectionf} that the original proof in loc.cit. fails, and in Appendix B that the proof of Tamarkin [T] fails as well; the failure in general is proven in Appendix A).\footnote{Here we mean the failure of the $\gl(V)$-equivariant formality. As well, we mean the formality of the Hochshild complex $\Hoch_\fin(S(V^*))$, introduced in Section \ref{section1.2} below, as a dg Lie algebra. Note that the proof in Appendix A requires only a weaker equivariance.} In this paper, we formulate and prove a statement closely related to the Kontsevich formality of Hochschild cochain on $S(V^*)$, which holds for an infinite-dimensional $V$ as well. 

In this paper, we restrict ourselves by considering the invinite-dimensional vector spaces $V$ over $\mathbb{C}$, which are graded:
\begin{equation}\label{eq17_new2}
V=\oplus_{i\ge 0}V_i,\ \ \dim V_i<\infty
\end{equation}
and such that  all graded components
 $V_i$ are {\it finite-dimensional}. 
 
 We assume the vector space $V$ to have the homological degree 0. The degree $i$ of all elements in $V_i$ is called the {\it auxiliary degree}.

For such a vector space $V$, we define the algebra of polynomial functions on $V$ as 
\begin{equation}\label{eq17_new1}
 S(V^*):=S(\oplus_{i\ge 0}V_i^*)
 \end{equation}
 The algebra $S(V^*)$ inherits from $V$ the grading, with $\deg V_a^*=-a$.

We start with defining suitable versions of
polyvector fields on $V$ and of the cohomological Hochschild complex
of $S(V^*)$. Then we prove in our setting  a direct analog of the
Hochschild-Kostant-Rosenberg theorem.

\subsection{The polyvector fields $T_\fin(V)$}\label{section1.1}
We define a suitable analogue of the Lie algebra of polyvector
fields $T_\poly(V)$, which we denote $T_\fin(V)$ . We want to allow infinite sums of monomials. Here is
the precise definition. 

A $k$-polyvector field of the grading (aka the auxiliary degree) $i$ is an element in
\begin{equation}\label{eq17_new3}
T_\fin(V)^{k,i}=\{\gamma\in \prod_{a_1,\dots,a_k\ge 0}\big(S(V^*)\otimes \Lambda^{a_1,\dots,a_k}V\big)|\ \deg_\aux\gamma=i\}
\end{equation}
where
$$
\Lambda^{a_1,\dots,a_k}V=V_{a_1}\wedge V_{a_2}\wedge\dots\wedge V_{a_k}
$$
Note that $\Lambda^{a_1,\dots,a_k}V$ is finite-dimensional, because all the graded components $V_i$ are. The grading $\deg_\aux$ is defined just below.

Note that  any element in $S(V^*)$ is a finite sum, by the definition \eqref{eq17_new1}.

The grading $\deg_\aux$ is defined for $f\otimes\lambda\in S(V^*)\otimes \Lambda^{a_1,\dots,a_k}V$ as
$$
\deg_\aux(f\otimes\lambda)=\deg f+a_1+\dots+a_k
$$
(note that the grading of an arbitary element in $\Lambda^{a_1,\dots,a_k}V$ is $a_1+\dots+a_k$, and the grading of an arbitary element in $V_a^*=-a$).

Then we define
\begin{equation}\label{eq17_new4}
T_\fin^k(V)=\bigoplus_iT_\fin^{k,i}(V)
\end{equation}

We denote
\begin{equation}\label{eq17_new5}
T_\fin^\udot(V)=\bigoplus_kT_\fin^k(V)[-k]
\end{equation}

\begin{example}\label{example1}
{\rm
Let $\frak{g}=\oplus_i\frak{g}_i$ be a graded Lie algebra, with all components $\frak{g}_i$ finite-dimensional (but $\frak{g}$ may be infinite-dimensional). Then the Kostant-Kirillov bivector of $\frak{g}$ is an element of $T_\fin^2(\frak{g}^*)$. Here $\frak{g}^*:=\oplus_i\frak{g}_i^*$. This Kostant-Kirillov bivector is a product of monomials, each of which has the grading (the auxilary degree) 0. }
\end{example}

Recall the Schouten-Nijenhuis bracket of polyvector fields on a (smooth manifold, smooth algebraic variety) $M$. It is a graded Lie bracket on $T_\poly(M)[1]$. The classical counterpart of the Schouten-Nijenhuis bracket we consider here is the one on the polynomial polyvector fields on a finite-dimensional vector space $W$, $T_\poly(W)$. For $\gamma_1,\gamma_2\in T_\poly(W)$, the Schouten-Nijenhuis bracket is
\begin{equation}\label{eq17_new1}
[\gamma_1,\gamma_2]=\gamma_1\circ \gamma_2-(-1)^{(p-1)(q-1)}\gamma_2\circ\gamma_1
\end{equation}
(where $\gamma_1\in T_\poly^p(W)$, $\gamma_2\in T_\poly^q(W)$).

The operation $\gamma_1\circ \gamma_2$ is defined as follows. Let $\gamma_1\in S^a(W^*)\otimes \Lambda^p(W), \gamma_2\in S^b(W^*)\otimes \Lambda^q(W)$. Then the operation $\gamma_1\circ \gamma_2$ is an operation
$$
-\circ-\colon \Big(S^a(W^*)\otimes \Lambda^p(W)\Big)\otimes\Big(S^b(W^*)\otimes\Lambda^q(W)\Big)\to S^{a+b-1}(W^*)\otimes\Lambda^{p+q-1}W
$$
There is the canonical $\gl(W)$-invariant $e^*\in (W\otimes W^*)^*$. Within the canonical isomorphism $W\otimes W^*\eqto \End(W)$, the operator $e^*$ is corresponded to the trace operator $\Tr\colon \End(W)\to\mathbb{C}$. The assumption that $\dim W<\infty$ is crucial for $e^*$ to exist.

The operator $e^*$ can be extended in a unique way to an $\gl(W)$-invariant operator $$e^*_{pb}\colon  S^b(W^*)\otimes \Lambda^p(W)\to S^{b-1}(W^*)\otimes \Lambda^{p-1}W  $$

Now $\gamma_1\circ \gamma_2$ is defined in two steps. At the first step, we just take the tensor product of all four factors
$$
S^a(W^*)\otimes \underset{e^*_{pb}}{\underbrace{\Lambda^p(W)\otimes S^b(W^*)}}\otimes\Lambda^q(W)\to S^a(W^*)\otimes \big(\Lambda^{p-1}W\otimes S^{b-1}(W^*)\big)\otimes \Lambda^qW
$$
and apply the operator $e^*_{pb}$ to the two factors in the middle. After that, we take the product of all factors, to consider it as an element of $S^\udot(W^*)\otimes\Lambda^\udot(W)$:
\begin{equation}
\gamma_1\circ\gamma_2=(f_1\otimes \ell_1)\circ (f_2\otimes \ell_2)=f_1\cdot e^*_{pb}(\ell_1,f_2)\cdot \ell_2
\end{equation}
where $f_i$ and $\ell_i$ are the factors in $S^\udot(W^*)$ and in $\Lambda^\udot(W)$, correspondingly.

This definition of $\gamma_1\circ \gamma_2$ is not valid when $W$ is an infinite-dimensional vector space. The reason is that the elementary invariant $e^*\colon W\otimes W^*\to\mathbb{C}$ is not defined when $\dim W=\infty$.

However, for the graded space $T_\fin^\udot(V)$ the operation $\gamma_1\circ \gamma_2$ is still well-defined, despite the corresponding map $\Tr\colon V^*\otimes V\to\mathbb{C}$ does not exist.  Indeed, it follows from  \eqref{eq17_new3} that the coefficient at each exterior algebra component $\Lambda^{a_1,\dots,a_k}V$ is a {\it (finite) polynomial} in $S(V^*)$. Let 
\begin{equation}\label{eq17_new6}
\gamma_1\in  \prod_{a_1,\dots,a_k}\underset{\gamma_1^{a_1,\dots,a_k}}{S(V^*)\otimes\Lambda^{a_1,\dots,a_k}V},\ \ \gamma_2\in\prod_{b_1,\dots,b_\ell}\underset{\gamma_2^{b_1,\dots,b_\ell}}{S(V^*)\otimes\Lambda^{b_1,\dots,b_\ell}V}
\end{equation}
Denote by $\gamma_1^{a_1,\dots,a_k}$, $\gamma_2^{b_1,\dots,b_\ell}$ the corresponding components of $\gamma_1$ and $\gamma_2$, correspondingly.

We claim that the following definition makes sense:
\begin{equation}\label{eq17_new7}
\gamma_1\circ\gamma_2=\prod_{\substack{{a_1,\dots,a_k}\\{b_1,\dots,b_\ell}}}\gamma_1^{a_1,\dots,a_k}\circ \gamma_2^{b_1,\dots,b_\ell}
\end{equation}
That is, we claim that the coefficient at {\it any fixed} component $\Lambda^{c_1,\dots,c_{k+\ell-1}}V$ is a {\it finite sum} of polynomials.

Assume that $\gamma_1^{a_1,\dots,a_k}\circ\gamma_2^{b_1,\dots,b_\ell}$ contributes to $S(V^*)\otimes \Lambda^{c_1,\dots,c_{k+\ell-1}}$.
Then one has $\{b_1,\dots,b_\ell\}\subset\{c_1,\dots,c_{k+\ell-1}\}$, the remaining $p-1$ of $\{c_1,\dots,c_{k+\ell-1}\}$ are in $\{a_1,\dots,a_k\}$. That is, (for any particular summand) there is exactly 1 index $i$ such that $a_i\not\in\{c_1,\dots,c_{k+\ell-1}\}$.
As the coefficient within $\Lambda^{b_1,\dots,b_\ell}V$ is a polinomial $f_{b_1,\dots,b_\ell}\in S(V^*)$, we get a condition on $\{a_1,\dots,a_k\}$: one has
\begin{equation}
\frac{\partial}{\partial x_{a_j}}f_{b_1,\dots,b_\ell}\ne 0\ \ \text{for some $j$ and $x_{a_j}\in V_{a_j}^*$}
\end{equation}
This condition depicts a finite set of possiblities for $a_j$, for fixed $\{b_1,\dots,b_\ell\}$ and $\{c_1,\dots,c_{k+\ell-1}\}$. As well, there is a finite number of choices of $\{b_1,\dots,b_\ell\}\subset
\{c_1,\dots, c_{k+\ell-1}\}$ (for a fixed $\{c_1,\dots,c_{k+\ell-1}\}$). It proves that the coefficients in  \eqref{eq17_new7} are finite sums, what makes sense of \eqref{eq17_new7}.

Then we define the Schouten-Nijenhuis bracket by \eqref{eq17_new1}. One easily sees that in this way one gets a graded Lie bracket on $T_\fin^\udot(V)[1]$. We have:
\begin{lemma}\label{lemma1}
Let $V$ be as in \eqref{eq17_new2}, and let $T_\fin^\udot(V)$ be as above.  Then the construction of the circle operation $\gamma_1\circ\gamma_2$, $\gamma_1,\gamma_2\in T_\fin^\udot(V)$ defines, via \eqref{eq17_new1}, a graded Lie bracket on $T_\fin(V)[1]$.
\end{lemma}
\qed

The graded Lie bracket on $T_\fin^\udot(V)[1]$, given in Lemma \ref{lemma1}, is called {\it the Schouten-Nijenhuis Lie bracket}.

\begin{example}\label{example2}
{\rm
Let us consider {\it vector fields} $v_1,v_2\in T^\udot_\fin(V)$, moreover, {\it linear} vector fields. That is, $$v_1,v_2\in \prod_a (\oplus_b V_b^*)\otimes V_a$$
In this case, we recover the product of the generalized Jacobian matrices.}
\end{example}

\begin{remark}\label{remark1}
{\rm
Note that all results of this Subsection would remain true if we defined $T_\fin^\udot(V)$ as
\begin{equation}
\overset{\sim}{T^\udot}_\fin(V)=\prod_{a_1,\dots,a_\ell}S(V^*)\otimes \Lambda^{a_1,\dots,a_\ell}V
\end{equation}
thus dropping the assumption on the auxilary degree in \eqref{eq17_new3} and \eqref{eq17_new5}. We adopt our previous definition, by the following reason. Later on, we define the corresponding ``finite'' version of the Hochschild cohomological complex $\Hoch^\udot_\fin(S(V^*))$, and prove the analogue of the Hochschild-Kostant-Rosenberg theorem. That is, one has 
$$
H^\udot(\Hoch^\udot_\fin(S(V^*)))=T^\udot_\fin(V)
$$
We do not know any definition of $\overset{\sim}{\Hoch^\udot}_\fin(S(V^*))$ such that one had
$$
H^\udot(\overset{\sim}{\Hoch^\udot}_\fin(S(V^*)))=\overset{\sim}{T^\udot}_\fin(V)
$$
}
\end{remark}

\subsection{The Hochschild complex $\Hoch^\mb_\fin(S(V^*))$}\label{section1.2}
Here we define a suitable version of the Hochschild cohomological complex of the algebra
$S(V^*)=S(\oplus_i V_i^*)$. We want its cohomology to be equal to $T^\udot_\fin(V)$.

Define
\begin{equation}
\begin{aligned}
\ &\Hoch^{k}_\fin(S(V^*))=\\
&\{\Psi\in \Hom(S(V^*)^{\otimes k},S(V^*))| \ \exists N(\Psi)\in \mathbb{Z} \text{  such that  }\\
&|\deg\Psi(f_1\otimes\dots\otimes f_k)-\sum_{i=1}^k\deg f_k|\le N(\Psi)\text{  for all  homogeneous  } f_1,\dots,f_k\}
\end{aligned}
\end{equation}
In other words, a cochain from $\Hoch^{k}_\fin(S(V^*))$ preserves the degree, up to a finite integral number, depending on the cochain.

One easily sees that $\Hoch^k_\fin(S(V^*))$ is a vector space. It inherits the grading, so that
$$
\Hoch^k_\fin(S(V^*))=\oplus_i\Hoch^{k,i}_\fin(S(V^*))
$$
where
$$
\begin{aligned}
\ &\Hoch^{k,i}(S(V^*))=\\
&\{\Psi\in\Hoch^k_\fin(S(V^*))|\ \deg\Psi(f_1\otimes\dots\otimes f_k)=\sum_{i=1}^k\deg f_k+i \text{  for all homogeneous  }f_1,\dots,f_k\}
\end{aligned}
$$

For any fixed $i$, 
$$
\Hoch^{\udot,i}_\fin(S(V^*))=\bigoplus_k\Hoch_\fin^{k,i}(S(V^*))[-k]
$$
becomes a complex with the Hochschild differential. Therefore, 
$$\Hoch^\udot_{\fin}(S(V^*))=\bigoplus_{k,i}\Hoch_\fin^{k,i}(S(V^*))[-k]$$ 
is a complex.

One easily sees that $\Hoch^\udot_{\fin}(S(V^*))[1]$ is a dg Lie algebra, with the Gerstenhaber Lie bracket.

The cohomology of $\Hoch^\udot_{\fin}(S(V^*))$ can be interpreted as a
derived functor, as follows.

Denote by $\EuScript{A}$ the category whose objects are {\it graded} $S(V^*)$-bimodules, and whose
morphisms are the grading preserving maps of them. Then $\EuScript{A}$ is an
abelian category. For an object $X\in \mathrm{Ob}(\EuScript{A})$,
denote by $X\langle j \rangle$ the object of $\EuScript{A}$ whose
inner grading is shifted by $j$: $(X\langle j\rangle)^i=X^{j+i}$.

We have the following lemma:
\begin{lemma}
The cohomology $H^\ell(\Hoch^\udot_{\fin}(S(V^*)))$ is equal to the vector space
$$H^\ell(\Hoch^\udot_{\fin,\tot}(S(V^*)))=\Ext^\ell_{\EuScript{A}}(S(V^*),\oplus_{j\in
\mathbb{Z}}S(V^*)\langle j\rangle)$$ 
\begin{proof}
The bar-resolution of $S(V^*)$ is clearly a projective resolution  in $\EuScript{A}$ of
the tautological bimodule $S(V^*)$. We compute the
$\Ext$'s functors by making use of this resolution. The complex
$\Hom_{\EuScript{A}}(\mathrm{Bar}^\mb(S(V^*)),\oplus_{k\in
\mathbb{Z}}S(V^*)\langle k\rangle)$ is exactly the complex
$\Hoch^\mb_{\fin}(S(V^*))$.
\end{proof}
\end{lemma}

\subsection{The Hochschild-Kostant-Rosenberg theorem}
Define the Hochschild-Kostant-Rosenberg map $\varphi_{HKR}\colon
T_\fin^\udot(V)\to \Hoch^\mb_\fin(V)$ as
\begin{equation}\label{1.3.1}
\varphi_{HKR}(\gamma)=\frac1{k!}\bigl\{f_1\otimes\dots\otimes
f_k\mapsto \gamma(df_1\wedge\dots\wedge df_k)\bigr\}
\end{equation}
for $\gamma\in T^k_\fin(V)$.

We have:
\begin{theorem}\label{theoremhkr}
The map $\varphi_{HKR}\colon T_\fin^\udot(V)\to\Hoch^\mb_\fin(S(V^*))$ is
a quasi-isomorphism of complexes.
\begin{proof}
Consider the following Koszul complex $K^\mb$:
\begin{equation}\label{1.3.2}
\dots\xrightarrow{d}K_3\xrightarrow{d}K_2\xrightarrow{d}K_1\xrightarrow{d}K_0\rightarrow
0
\end{equation}
where
\begin{equation}\label{1.3.3}
K_k=\bigoplus_{a_1,\dots,a_k} S(V^*\oplus V^*)\otimes\Lambda^{a_1,\dots,a_k}V^*
\end{equation}
with the
differential 
\begin{equation}\label{1.3.4}
d((\xi_{i_1}\wedge\dots\wedge \xi_{i_k})\otimes
f)=\sum_{j=1}^k(-1)^{j-1}(\xi_{i_1}\wedge\dots\wedge\hat{\xi}_{i_j}\wedge\dots\wedge
\xi_{i_k})\otimes((x_j-y_j)f)
\end{equation}
where $\{x_i\}$ is a basis in $V^*$ compatible with the decomposition
$V=\oplus_i V_i$, $\{y_i\}$ is the same basis in the second copy of
$V^*$, and $\{\xi_i\}$ is the corresponded basis in $V[1]$.

This Koszul complex is clearly a resolution of the tautological
$S(V^*)$-bimodule $S(V^*)$ by free bimodules.  It is as well a resolution in
the category $\EuScript{A}$, because the differential $d$ preserves
the auxilary grading. We compute:
\begin{equation}\label{eq174}
\begin{aligned}
\ &\Hom_{\EuScript{A}}\Big(K^\ell,\bigoplus_{j\in \mathbb{Z}}S(V^*)\langle
j\rangle\Big)=\\&\Hom_{\mathbb{C}}\Big(\bigoplus_{a_1,\dots,a_\ell}\Lambda^{a_1,\dots,a_\ell}V^*, \bigoplus_{j\in \mathbb{Z}}S(V^*)\langle
j\rangle\Big)=\prod_{a_1,\dots,a_\ell}\Hom_{\mathbb{C}}\Big(\Lambda^{a_1,\dots,a_\ell}V^*,\bigoplus_{j\in \mathbb{Z}}S(V^*)\langle
j\rangle\Big)=\\
&\prod_{a_1,\dots,a_\ell}\bigoplus_j[S(V^*)\otimes \Lambda^{a_1,\dots,a_\ell}V]_j
\end{aligned}
\end{equation}
where the subscript $j$ denotes the elements of the auxilary grading $j$. In the last equality we essentially use that all $\Lambda^{a_1,\dots,a_\ell}V^*$ are finite-dimensional.

The rightmost term in \eqref{eq174} is exactly $T_\fin^\ell(V)$.

One easily checks that the induced differential vanishes. It completes the computation of the cohomology $H^\udot(\Hoch_\fin^\udot(S(V^*)))$.

It remains to note that the image of the
Hochschild-Kostant-Rosenberg map $\varphi_{HKR}$ coincides with the cohomology
classes in $\Hoch^\mb_\fin(V)$ produced by the Koszul resolution.
\end{proof}
\end{theorem}

\subsection{The polydifferential operators, associated with graphs}\label{2017graphs}
Here we recall the construction from [K97], assigning a Hochschild cochain of the algebra $S(W^*)$, to several polyvector fields on $W$, and to a combinatorial datum given by a {\it Kontsevich admissible graph}. We refer the reader to [K97, Section 6.1] for more detail.

Throughout this Subsection, $W$ denotes a finite-dimensional vector space over $\mathbb{C}$.

The goal is to construct the most general Hochschild cochains, associated with an ordered sequence of polyvector fields $\gamma_1,\dots,\gamma_n$ on $W$, such that the construction is  $\gl(W)$-equivariant. 

The idea is to subsequently apply the elementary invariant operator $e^*\colon (W\otimes W^*)\to\mathbb{C}$, associated to all edges of the graph. In fact, the construction is a generalization of the construction of $\gamma_1\circ\gamma_2$ (see Section \ref{section1.1}), which is corresponded to the graph with two vertices ``of the first type'' (see below), and a single oriented edge, see Figure \ref{figurenew17.0}.

\sevafigc{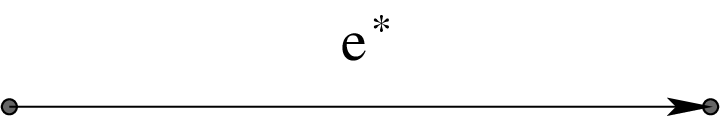}{50mm}{0}{To each oriented edge, is associated the elementary invariant $e^*\in (W\otimes W^*)^*\sim W^*\otimes W$\label{figurenew17.0}}

For a general graph $\Gamma$, we use the notation $V_\Gamma$ for the set of vertices of $\Gamma$, and $E_\Gamma$ for the set of its edges. 
\begin{defn}\label{kontgraphs}{\rm
A Kontsevich admissible graph $\Gamma$ is an oriented graph with two types of labelled vertices, the vertices of the first type, labelled $\{1,\dots,n\}$, and the vertices of the second type, labelled $\{\bar{1},\dots,\bar{m}\}$, such that
\begin{itemize}
\item[(i)] $\sharp E_\Gamma= 2n+m-2\ge 0$, $V_\Gamma=\{1,\dots,n\}\sqcup\{\bar{1},\dots,\bar{m}\}$, 
\item[(ii)] every edge $(v_1,v_2)\in E_\Gamma$ starts at a vertex of the first type, $v_1\in\{1,\dots,n\}$,
\item[(iii)] there are no simple loops (aka tadpoles), that is, edges of the form $(v,v)$,
\end{itemize}
}
\end{defn}
\begin{remark}\label{remark171new}{\rm
Our definition is slightly different from the original one, as we do not consider an ordering of all sets $\Star(v)$, $v$ a vertex of the first type, as a part of the data. The reason is that, as soon as the sets of  all vertices of the first type and of all vertices of the second type are ordered, there is an induced ordering of the edges.
}
\end{remark}
For a Kontsevich admissible graph $\Gamma$, we denote by $V_\Gamma^I$ the set of vertices of the first type, and by $V_\Gamma^{II}$ the set of vertices of the second type, so that $V_\Gamma=V_\Gamma^I\sqcup V_\Gamma^{II}$.

Let us recall the construction, associating to a Kontsevich admissible graph $\Gamma$ with $n$ vertices of the first type and $m$ vertices of the second type, $n$ homogeneous polyvector fields $\gamma_1,\dots,\gamma_n$ on $W$, and $m$ functions $f_1,\dots,f_m\in S(W^*)$, a {\it polyvector field} $\mathcal{U}_\Gamma(\gamma_1,\dots,\gamma_n;\ f_1,\dots,f_m)$ of the (cohomological) degree 
\begin{equation}
\deg \mathcal{U}_\Gamma(\gamma_1,\dots,\gamma_n;\ f_1,\dots,f_m)=\sum_{i=1}^n\deg \gamma_i -2n-m+2
\end{equation}
When
\begin{equation}\label{eq17_new10}
\sum_{i=1}^n\deg\gamma_i=2n+m-2
\end{equation}
the polyvector field $\mathcal{U}_\Gamma(\gamma_1,\dots,\gamma_n;\ f_1,\dots,f_m)$ has (cohomological) degree 0, that is, is a function. In this case, we get a {\it Hochschild cochain} $\mathcal{U}_\Gamma(\gamma_1,\dots,\gamma_n)\in \Hom_\mathbb{C}(S(W^*)^{\otimes m},S(W^*))$.

A polynomial polyvector field on $W$ is an element of the graded commutative algebra $T_\poly(W)=S(W^*\oplus W[-1])$.
Let $\Gamma$ be a Kontsevich admissible graph, with $\sharp V_\Gamma^I=n$, $\sharp V_\Gamma^{II}=m$.
Consider the associative graded commutative algebra
\begin{equation}
A_\Gamma=\prod_{v\in V_\Gamma^I}(S(W^*\oplus W[-1]))_v\ \otimes\ \prod_{v\in V_\Gamma^{II}}(S(W^*))_v
\end{equation}
Although the construction of $\Gamma\rightsquigarrow \mathcal{U}_\Gamma$ does not depend on the choice of basis in $W$, we choose one a write the construction ``in coordinates'', as it makes it more readable.

Let $\{x_i\}$, $i=1,\dots, N=\dim W$ be a basis in $W^*$, let $\{\xi_i^\prime\}$, $i=1,\dots,N$ be the dual basis in $W$, and let $\{\xi_i\}$, $i=1,\dots,N$ be the corresponding basis in $W[1]$.

One assigns to the elementary invariant $e^*\in (W\otimes W^*)^*$ the operator 
\begin{equation}
e^*=\sum_{i=1}^N\frac{\partial}{\partial x_i}\otimes \frac{\partial}{\partial \xi_i}
\end{equation}
which does not depend on the choice of the basis $\{x_i\}$. 

Let $(v_1,v_2)$ be an oriented edge of $\Gamma$. We assign to it the operator acting in $A_\Gamma$:
\begin{equation}
e_{(v_1,v_2)}^*=\sum_{i=1}^N(\frac{\partial}{\partial x_i})_{v_2}\otimes (\frac{\partial}{\partial \xi_i})_{v_1}
\end{equation}
where the sub-indices $v_1$ and $v_2$ indicate the factors in $A_\Gamma$ on which the corresponding operators act.

The operators $e^*_{(v_1,v_2)}$, acting on $A_\Gamma$, commute up to a sign, for different edges. The labellings which is a part of the definition of an admissible graph, fix an order on all edges. In the formula below this order is assumed:
\begin{equation}
D_\Gamma=\prod_{(v_1,v_2)\in E_\Gamma}e^*_{(v_1,v_2)}\colon A_\Gamma\to A_\Gamma
\end{equation}
Take homogeneous polyvector fields $\gamma_1,\dots,\gamma_n$, and functions $f_1,\dots,f_m$. The ordering of the sets $V_\Gamma^I$ and $V_\Gamma^{II}$ fixes an ordering of them. In the formula below we assume this ordering:
\begin{equation}
\mathcal{U}_\Gamma(\gamma_1,\dots,\gamma_n)(f_1,\dots,f_m)=\Delta^*\Big(D_\Gamma(\gamma_1\otimes\dots\otimes\gamma_n\otimes f_1\otimes\dots\otimes f_m)\Big)
\end{equation}
where $\Delta^*(-)$ is the ``restriction of the function to the diagonal'', which is, by definition, the product of the components:
\begin{equation}
\Delta^*\colon A_\Gamma\to S(W^*\oplus W[-1])
\end{equation}

In general, $\mathcal{U}_\Gamma(\gamma_1,\dots,\gamma_n)(f_1,\dots,f_m)\in T_\poly(W)$ is a polyvector field. When \eqref{eq17_new10} holds, it is a function. In this case, $\mathcal{U}_\Gamma(\gamma_1,\dots,\gamma_n)$ is a Hochschild cochain of $S(W^*)$.

\subsection{Towards the formality for an infinite-dimensional space}\label{sectionf}
Here we explain why the Kontsevich solution of the formality theorem fails for the the algebra of polynomials on an infinite-dimensional vector space $S(V^*)$, in the sense of Section \ref{section1.0}. 
One of the main parts of the construction in [K97] is an assignment associating to each Kontsevich admissible graph $\Gamma$, and $n=\sharp V_\Gamma^I$ of homogeneous polyvector fields $\gamma_1,\dots,\gamma_n$ on $W$, with \eqref{eq17_new10}, the Hochschild cochain $\mathcal{U}_\Gamma(\gamma_1,\dots,\gamma_n)\in \Hoch^\udot(S(W^*))$, see Section \ref{2017graphs}.

Consider the infinite-dimensional framework of Sections \ref{section1.0}, \ref{section1.1}. 

We claim that, for a
graph $\Gamma$ containing a fragment which is {\it an oriented cycle between the vertices
of the first type}, the corresponding polydifferential operator
$\U_\Gamma$ is in general ill-defined.
An informal explanation is that the polydifferential operator corresponding to an oriented cycle ``looks like a trace operator on an infinite-dimensional vector space'', which
is ill-defined when $\dim V=\infty$. See Example \ref{example5} below, for an explicit computation.

In fact, the graphs with oriented cycles between the vertices of the first type are the only ``bad'' graphs: 
\begin{lemma}\label{lemma1714}
Let $\Gamma$ be a Kontsevich admissible graph in the sense of [K97]
with $n$ vertices of the first type and $m$ vertices of the second
type. Suppose that the graph $\Gamma$ does not contain any oriented cycle
between the vertices of the first type, as its sub-graph. Let $V$ be as above, and let
$\gamma_1,\dots,\gamma_n\in T_\fin(V)$. Then the Kontsevich
polydifferential operator
$\U_\Gamma(\gamma_1\wedge\dots\wedge\gamma_n)$ is well-defined
as an element of $\Hoch^\mb_\fin(S(V^*))$.
\end{lemma}

{\it Proof of Lemma:}

In our definitions \eqref{eq17_new2} and \eqref{eq17_new1}, an element of $S(V^*)$ is a {\it finite sum}. In contrast, an element of $T_\fin(V)$ is an {\it infinite product}, see \eqref{eq17_new3}.

Let $\Gamma$ be a Kontsevich admissible graph with oriented cycles between the vertices of the first type.
Then there is a vertex of the first type, say $v_0$, such that all edges outgoing from $v_0$ have as their targets vertices of the second type. 

Let $\gamma_1,\dots,\gamma_n$ and $f_1,\dots,f_m$ be fixed. Consider all operators $e^*_t=(\frac{\partial}{\partial \xi_i})_{v_0}\otimes (\frac{\partial}{\partial x_i})_{v_t}$, associated with all edges $t$ outgoing from $v_0$. As the functions $f_1,\dots,f_m$ are finite sums, the operator 
$$
\prod_{t=(v_0,v_t)}e^*_t(\gamma_{v_0}\otimes f_1\otimes\dots\otimes f_m)
$$
is non-zero only for a finitely many components $\gamma_{v_0}^{(a_1,\dots,a_\ell)}\in S(V^*)\otimes\Lambda^{(a_1,\dots,a_\ell)}V$.
Namelly, the numbers $a_1,\dots,a_n$ should obey the condition:
\begin{equation}\label{eq17_new15}
\frac{\partial}{\partial x_{a_i}}(f_1\cdot\dots\cdot f_m)\ne 0\text{    for $i=1,\dots,\ell$}
\end{equation}
For each such component, denote
$$
g(\gamma_{v_0})^{(a_1,\dots,a_\ell)}=\gamma_{v_0}(dx_{a_1}\wedge\dots\wedge dx_{a_\ell})
$$
Then we remove the vertex $v_0$ from $\Gamma$, as well as all incoming to and outgoing from $v_0$ edges.
Denote the obtained graph by $\Gamma^{(1)}$. The graph $\Gamma^{(1)}$ fulfils the assumptions of the Lemma as well.
We find a vertex $v_1$ of the first type of $\Gamma^{(1)}$ such that any edge outgoing from $v_1$ targets at a vertex of the second type. The components $\gamma_{v_1}^{(b_1,\dots,b_s)}\in S(V^*)\otimes \Lambda^{(b_1,\dots,b_s)}$ which contribute to $\mathcal{U}_\Gamma(\gamma_1,\dots,\gamma_n)(f_1,\dots,f_m)$ by a non-zero summand are those for which
\begin{equation}\label{eq17_new16}
\frac{\partial}{\partial x_{b_i}}(g(\gamma_{v_0})^{(a_1,\dots,a_\ell)}\cdot f_1\cdot\dots\cdot f_m)\ne 0\text{  for all $i=1\dots s$ and for some $(a_1,\dots,a_\ell)$ obeying \eqref{eq17_new15}}
\end{equation}
The set of all possible $(a_1,\dots,a_\ell)$ and $(b_1,\dots,b_s)$ which contribute by a non-zero summand is thus finite. Then we remove the vertex $v_1$ from $\Gamma^{(1)}$, with all its incoming and outgoing edges, and so on.

After successive iteration of this procedure, we get only vertices of the second type, and a finite number of summands which do finally contribute to $\mathcal{U}_\Gamma(\gamma_1,\dots,\gamma_n)(f_1,\dots,f_m)$.

\qed

\begin{example}\label{example5}{\rm
Here we consider an example of a graph $\Gamma$ with oriented cycles between the vertices of the first type, and of an infinite-dimensional space $V$, such that the Hochschild cochain $\mathcal{U}_\Gamma(\gamma_1,\dots,\gamma_n)$ is ill-defined for some $\gamma_1,\dots,\gamma_n\in T_\fin(V)$.

Consider the graph with $\sharp V_\Gamma^I=2, \sharp V_\Gamma^{II}=2$ shown in Figure \ref{figurenew17.1}. (The labelling is not essential for this Example).

\sevafigc{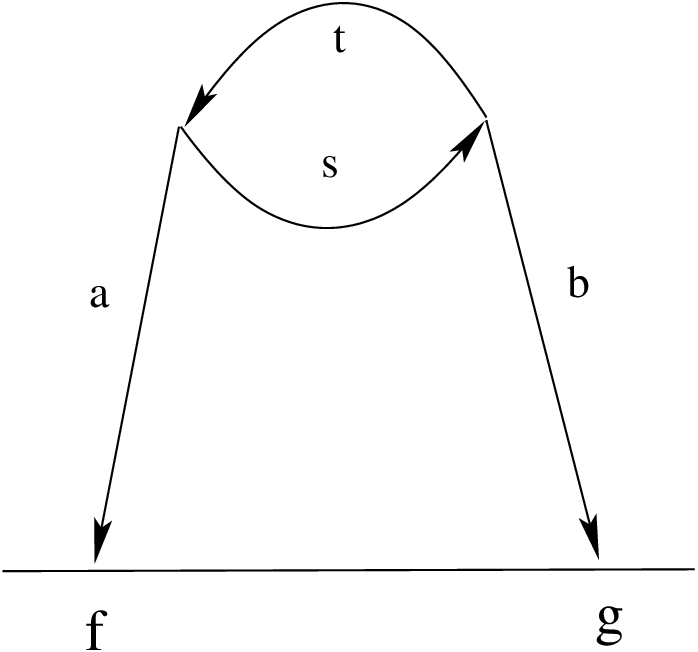}{50mm}{0}{A graph $\Gamma$ with an oriented cycle between the type I vertices, which leads to divergent operator $\mathcal{U}_\Gamma$\label{figurenew17.1}}

Consider the vector space $V$ as in Section \ref{section1.0} such that
$$
V=\oplus_{a\ge 0}V_a,\ \ \dim V_a=1\text{  for all $a$}
$$
We put in the vertices of the second type linear functions, for simplicity. Say, $f=x_a in V_a^*$, $g=x_b \in V_b^*$. We put in the vertices of the first type (the same) quadratic polyvector field $\alpha$ in $T_\fin(V)$,
$$
\alpha=\sum_{0\le i<j\le \infty} \sum_{s+t=i+j}x_sx_t\partial_i\wedge\partial_j
$$
where $x_a$ is a non-zero vector in $V_a^*$ (thus $\{x_a\}_{a\ge 0}$ is a basis in $V^*=\oplus V_a^*$), and $\{\partial_a\}_{a\ge 0}$ is the corresponding dual basis in $V$. 

We have:
$$
\begin{aligned}
\ &
\mathcal{U}_\Gamma(\alpha,\alpha)(x_a,x_b)=
\sum_{\substack{{s>a,t>b}\\
{s+a-t\ge 0, t+b-s\ge 0}}}x_{s+a-t}x_{t+b-s}-\sum_{\substack{{s<a,t>b}\\
{s+a-t\ge 0, t+b-s\ge 0}}}x_{s+a-t}x_{t+b-s}-\\
&\sum_{\substack{{s>a,t<b}\\
{s+a-t\ge 0, t+b-s\ge 0}}}x_{s+a-t}x_{t+b-s}+\sum_{\substack{{s<a,t<b}\\
{s+a-t\ge 0, t+b-s\ge 0}}}x_{s+a-t}x_{t+b-s}
\end{aligned}
$$
The condition $-a\le s-t\le b$ (common for all four summands) shows that the three last summands are in fact finite. 

On the other hand, the first summand is infinite, for any $a,b\ge 0$. 

}
\end{example}

\hspace{3pt}

The Lemma above shows that the only graphs for which the operators $\mathcal{U}_\Gamma$ are ill-defined for an infinite-dimensional $V$ (in the sense of Section \ref{section1.0}) are the graphs with oriented cycles between the vertices of the first type. Therefore, in order to extend the Kontsevich proof to the case of an infinite-dimensional $V$, we should ``exclude'' all these graphs from the consideration.
In this paper, we
suggest a way how to do that, by providing a new propagator
1-form, essentially different from those used by Kontsevich. Within this new propagator, one can mimic the most of Kontsevich's arguments in his proof of formality [K97], with one notable exception. Namely,  Lemma in loc.cit.,
Section 6.6 (saying that when $n\ge 3$ points of the first type approach each other far from the real line, the corresponding integral vanishes) fails for our new propagator. It implies that one gets a new
$L_\infty$ structure on $T_\fin(V)$, called here  {\it exotic} and denoted by $T_\fin^\mathcal{L}(V)$, and gets a construction  of an $L_\infty$ quasi-isomorphism from this exotic $L_\infty$ algebra $T_\fin^\mathcal{L}(V)$ to the Hochschild complex $\Hoch^\udot_\fin(S(V^*))$. The second Taylor
component of $T_\fin^\mathcal{L}(V)$ is given by the Schouten-Nijenhuis bracket, but there are higher
non-vanishing Taylor components as well. In fact, we show that all odd Taylor components vanish, and all even Taylor component are non-zero.

First of all, we provide the construction of this exotic $L_\infty$ structure. After that, we construct the $L_\infty$ morphism.

\section{The  $L_\infty$ algebra $T_\fin^\mathcal{L}(V)$}
\subsection{Configuration spaces $C_{n,\Gamma}$}
Recall that, for an oriented graph $\Gamma$, we denote by $V(\Gamma)$ the set of
its vertices, and by $E(\Gamma)$ the set of its edges. For a vertex
$v\in V(\Gamma)$, denote by $Star(v)$ the set of edges outgoing 
from the vertex $v$, and by $In(v)$ the set edges incoming to
the vertex $v$.

In this paper, we deal with three different types of admissible graphs. These are the {\it Kontsevich admissible graphs} (we recalled their Definition in Section \ref{2017graphs}), the {\it normalized plane admissible graphs}, introduced just below, and the {\it normalized Kontsevich admissible graphs}, see Section \ref{sectionnewprop}. 

\begin{defn}\label{defn1721}{\rm
A {\it normalized (plane) admissible graph} $\Gamma$ is
an oriented graph without any oriented loops, whose vertices are labelled by $\{1,2,\dots,\sharp V(\Gamma)\}$. 
Denote by $n(v)$ be the label of
vertex $v$, $n(v)\in \{1,\dots \sharp V(\Gamma)\}$. The labelling is subject to
the following condition: 
\begin{equation}\label{Knew}
n(v_1)<n(v_2)\text{ if there is an edge $\overrightarrow{v_1v_2}$}
\end{equation}
We do not assume a strong admissible graph to be connected. 
}
\end{defn}
We denote the set of strong admissible graphs $\Gamma$ with $n$ vertices
and $E$ edges by $G(n,E)$.

Note that, unlike for the case of Kontsevich admissible graphs, the vertices of a strong admissible graph are all of the same type.
We imagine a strong admissible graph as placed over the plane $\mathbb{R}^2$ (while a Kontsevich admissible graph should be imagined as placed over the closed upper-half plane $\mathbb{R}^2_+$, such that the vertices of the first type are placed over the open upper half-plane, and the vertices of the second type are placed over the boundary line).

Any graph without loops admits at least one labelling obeying
(\ref{Knew}), but it may admit several different such
labelings, see Figure \ref{figure1}.

\sevafigc{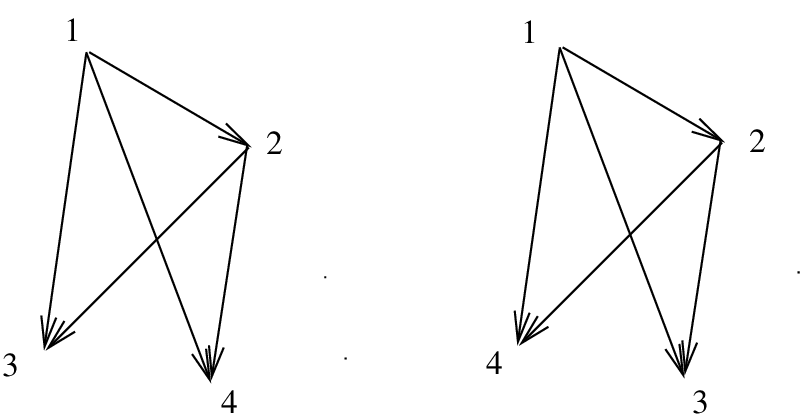}{80mm}{0}{The two different admissible labelings of a graph
with 4 vertices\label{figure1}}

Let $\Gamma\in G(n,E)$. Define the configuration space
$C_{n,\Gamma}$, as follows:
\begin{equation}\label{k2.1.1}
\begin{aligned}
\ &C_{n,\Gamma}=\{z_1,\dots,z_n,\, z_i\ne z_j\,\text{for $i\ne j$},\\
&\text{and $\mathrm{Im}(z_j-z_i)<0$ when $\overrightarrow{(i,j)}$ is
an edge of $\Gamma$}\}/G^3
\end{aligned}
\end{equation}
where $G^3=\{z\mapsto az+c,\, a\in\mathbb{R}_{>0},c\in \mathbb{C}\}$
is a 3-dimensional group of transformations.

In particular, if $E=0$ we recognize the Kontsevich's configuration
space $C_n$ from [K97]. Dimension $\dim C_{n,\Gamma}$ is equal to
$2n-3$ and does not depend on $\Gamma$.

Now suppose that $\sharp E(\Gamma)=2n-3$. If
$e=\overrightarrow{z_iz_j}$ is an edge of $\Gamma$, we associate
with it the differential 1-form
\begin{equation}\label{k2.1.2}
\phi_e=d\mathrm{Arg}(z_j-z_i)=\frac1{2i}d\mathrm{Log}\frac{z_j-z_i}{\overline{z}_j-\overline{z}_i}
\end{equation}
A labelling of the vertices of $\Gamma$ yields a 
labelling of the edges. To define it, one firstly takes the edges outgoing from
the vertex 1, in the order fixed by the labelling of their targets, then we take the
edges outgoing from the vertex 2, and so on.

Define
\begin{equation}\label{k2.1.3}
W_\Gamma=\frac1{\pi^{2n-3}}\int_{C_{n,\Gamma}}\bigwedge_{e\in
E(\Gamma)}\phi_e
\end{equation}
Here in the wedge product we use the order of the 1-forms $\phi_e$
corresponded to the ordering of the edges of $\Gamma$, described
just above.

\comment
\begin{remark}
In [K97], Section 6.2, M.Kontsevich considers the combinatorial factor
$\prod_{v\in V(\Gamma)}\frac1{\sharp Star(v)!}$ for an analogous
formula. The reason why we omit this factor here is the following.
In [K97] in the definition of an admissible graph not only the
vertices are labeled, but also the sets $Star(v)$ for all vertices
$v$ are ordered. It is not necessary because as soon as the vertices
are enumerated, there appears a canonical ordering, by the number of
the end-vertex. Therefore, we omit here and thereafter this
ingredient from the data of an admissible graph. As a consequence,
we have not the mentioned above combinatorial factor.
\end{remark}
\endcomment

We firstly show that $W_\Gamma=0$ for an odd $n$.
\begin{lemma}\label{lemma1721}
For any admissible graph with an odd number $n$ of vertices, the integral $W_\Gamma=0$ .
\begin{proof}
Map a vertex $p$ to the point $0+0i$, using the action of the group $G^3$, and
let another point $q$ move along the unit circle around $p$. (These two are the only vertices with the restricted degrees of freedom, in this way we get rid of the action of $G^3$). Draw
the vertical line $\ell$ through $p$ and consider the symmetry
$\sigma$ with respect to $\ell$. One has the following general
formula:
\begin{equation}\label{k2.1.4}
\int_c\sigma^*\omega=\int_{\sigma_*c}\omega
\end{equation}
where $c$ is an {\it oriented} chain.

In our case $c=C_{n,\Gamma}$ and $\omega=\bigwedge_{e\in
E(\Gamma)}\phi_e$. We have: $\sigma_*\phi_e=-\phi_e$ (there are
$2n-3$ edges), and at each ``movable'' point (there are $n-1$ of such
points) the orientation changes to the opposite. Finally, if $I$ is the
integral, we have from (\ref{k2.1.4}): $(-1)^{2n-3}I=(-1)^{n-1}I$
which implies that $I=0$ for an odd $n$.
\end{proof}
\end{lemma}

It is clear that when $\Gamma$ is the graph with two vertices and
a single edge, $W_\Gamma=1$. For an admissible graph $\Gamma$ with $n$ vertices to have a non-zero integral $W_\Gamma$, it should have $2n-3=\dim C_{n,\Gamma}$ edges. Consider several examples.

\begin{example}\label{example1721}{\rm  Consider the graph $\Gamma$ shown in Figure \ref{figure2}.\footnote{An example similar to this one was considered by M.Kontsevich, in his email to the author regarding the first archive version of the paper} 
\sevafigc{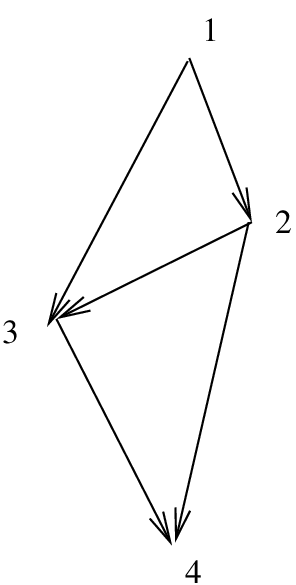}{30mm}{0}{An admissible graph $\Gamma$ with
$n=4$ and nonzero $W_\Gamma$\label{figure2}} Let us compute $W_\Gamma$ for this
graph. Fix the vertex 2 to the point $0+0\cdot i$ by the action of
group $G^3$, and let the vertex 3 move along the unit lower
half-circle around 2. Let $x$ be the angle of the arrow
$\overrightarrow{(2,3)}$, $-\pi\le x\le 0$. Then we can integrate
over positions of the vertices 1 and 4 separately. Let us firstly
integrate (for a fixed $x$) over 1, denote it $z$.

We need to compute the integral
\begin{equation}\label{ex1.1}
I_1=\int_{\mathrm{Im}z\ge 0}d\mathrm{Arg}(-z)\wedge
d\mathrm{Arg}(-z+\exp(ix))
\end{equation}
We use the Stokes formula:
\begin{equation}\label{ex1.2}
I_1=\int_{\partial(\mathrm{Im}z\ge
0)}\mathrm{Arg}(-z)d\mathrm{Arg}(-z+\exp(ix))
\end{equation}
The Stokes formula is applied to the domain $D=\{\mathrm{Im}z\ge
0,\,\varepsilon\le |z|\le R\}$ where $\varepsilon\to 0$ and
$R\to\infty$. There are two boundary strata which contribute to the
integral:
\begin{equation}\label{ex1.3}
I_1=-\lim_{\varepsilon\to 0}\int_{|z|=\varepsilon,\,\mathrm{Im}z\ge
0}\mathrm{Arg}(-z)d\mathrm{Arg}(-z+\exp(ix))+\lim_{R\to\infty}\int_{|z|=R,\,\mathrm{Im}z\ge
0}\mathrm{Arg}(-z)d\mathrm{Arg}(-z+\exp(ix))
\end{equation}
The first integral in the limit $\varepsilon\to 0$ is equal to $\pi
x$, and the second one in the limit $R\to\infty$ is equal to
$\frac12\pi^2$. The total answer is
\begin{equation}\label{ex1.4}
I_1=-\pi x+\frac12\pi^2
\end{equation}
Analogously we compute the integral $I_2$ over position of the
vertex 4. We get the same answer:
\begin{equation}\label{ex1.5}
I_2=-\pi x+\frac12\pi^2
\end{equation}
Finally,
\begin{equation}\label{ex1.6}
W_\Gamma=\frac1{\pi^5}\int_{-\pi}^0I_1(x)I_2(x)dx=\frac{13}{12}
\end{equation}
}
\end{example}

\begin{example}\label{example1722}{\rm More generally, consider the graph $\Gamma$ shown
in Figure \ref{figure3}. An analogous computation shows that
\begin{equation}\label{ex2.1}
\begin{aligned}
\ &W_\Gamma=\frac1{\pi^{2m+2n+1}}\int_{-\pi}^0(\pi
x-\frac12\pi^2)^{m+n}dx=\\
&\frac1{\pi^{2m+2n+1}}\cdot \Bigl(\frac1{\pi}\frac1{m+n+1}(\pi
x-\frac12\pi^2)^{m+n+1}\Bigr)|_{-\pi}^0=\\
&\begin{cases}
(-1)^{m+n}\frac{3^{m+n+1}-1}{(m+n+1)2^{m+n+1}}&\text{if  }m+n\text{  is even}\\
0&\text{if  }m+n\text{  is odd}
\end{cases}
\end{aligned}
\end{equation}
(The additional boundary strata, coming from the components when some $\ge 2$ of the $m$ upper points, or some $\ge 2$ of the $n$ lower points, approach each other, clearly do not contribute to the integral).
It follows that $W_\Gamma\ne 0$ for any $m,n$ such that $m+n$ is even.
\sevafigc{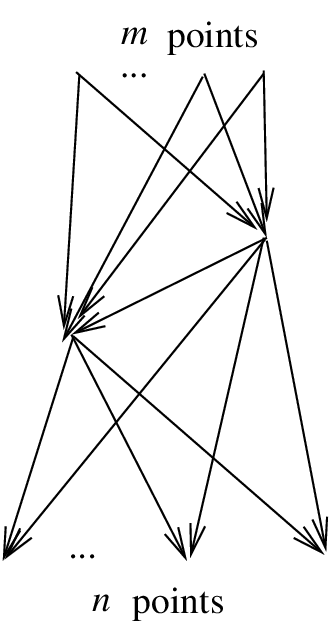}{35mm}{0}{The graph $\Gamma$ from Example 2\label{figure3}}
}
\end{example}

\subsection{The exotic $L_\infty$ structure}
We are going to define polylinear operators
$$
\mathcal{L}_n\colon \Lambda^nT_\fin(V)\to T_\fin(V)[2-n], \ n\ge 2
$$
which later are proven to be the Taylor components of an $L_\infty$ structure.

Let $\gamma_1,\dots,\gamma_n$ are homogeneous polyvector fields. We
are going to define the value
$\mathcal{L}_n(\gamma_1\wedge\dots\wedge\gamma_n)$.

First of all, we define with an admissible graph $\Gamma$ with $n$
vertices and for an ordered set $\gamma_1,\dots,\gamma_n$ of polyvector
fields a polyvector field
$\mathcal{L}_\Gamma(\gamma_1\otimes\dots\otimes\gamma_n)$ which is
$\bigl(\sum_{i=1}^n\deg\gamma_i+n-\sharp E(\Gamma)\bigr)$-vector
field. (Here we denote by $\deg\gamma$ the Lie algebra degree, that
is, if $V$ is concentrated in the cohomological degree zero,
$\deg\gamma=k-1$ for a $k$-vector field $\gamma$).

If the target vector space $V$ were finite-dimensional, the
polyvector field
$\mathcal{L}_\Gamma(\gamma_1\otimes\dots\otimes\gamma_n)$ would be the sum
\begin{equation}\label{k2.2.1}
\mathcal{L}_\Gamma(\gamma_1\otimes\dots\otimes\gamma_n)=\sum_{I\colon
E(\Gamma)\to\{1,2,\dots,\dim
V\}}\mathcal{L}_\Gamma^I(\gamma_1\otimes\dots\otimes\gamma_n)
\end{equation}
The polyvector field
$\mathcal{L}_\Gamma^I(\gamma_1\otimes\dots\otimes\gamma_n)$ is the
product over the vertices of $\Gamma$:
\begin{equation}\label{k2.2.2}
\mathcal{L}_\Gamma^I(\gamma_1\otimes\dots\otimes\gamma_n)=\bigwedge_{v\in
V(\Gamma)}\Psi_v^I
\end{equation}
where the vertices are taken in the order corresponded to the labelling of the
vertices. Each $\Psi_v$ is defined as
\begin{equation}\label{k2.2.3}
\Psi_v^I=\left(\prod_{e\in In(v)}\frac{\partial}{\partial
x_{I(e)}}\right)\bigl(\gamma_{n(v)},\wedge_{e\in Star(v)}dx_{I(e)}\bigr)
\end{equation}
where $n(v)$ is the label of the vertex $v$, and the order in the
wedge-product is fixed by the ordering of the set $Star(v)$, as
above (see Remark \ref{remark171new}). It completes the definition of
$\mathcal{L}_\Gamma(\gamma_1\otimes\dots\otimes\gamma_n)$ for a
finite-dimensional vector space $V$.

Now suppose that $V=\bigoplus_{i\in \mathbb{Z}_{\ge 0}}V_i$ where
all graded components $V_i$ are finite-dimensional, as in Section \ref{section1.0}.
We leave to the reader the following lemma:
\begin{lemma}
Let
$\gamma_i\in T_\fin(V)$, and a graph $\Gamma$ does not contain
any oriented cycles (this assumption is fulfilled if $\Gamma$ is an admissible graph).
Then the polyvector field
$\mathcal{L}_\Gamma(\gamma_1\otimes\dots\otimes\gamma_n)$ is
well-defined for the case of infinite-dimensional $V$, as above.
\end{lemma}
\qed

Let
$\gamma_1,\dots,\gamma_n\in T_\fin(V)$. Define the polyvector field in $T_\fin(V)$
\begin{equation}\label{k2.2.4}
\mathcal{L}_n(\gamma_1\wedge\dots\wedge\gamma_n)=\mathrm{Alt}_{\gamma_1,\dots,\gamma_n}\Bigl(\sum_{\Gamma\in
G_{n,2n-3}}W_\Gamma\cdot\mathcal{L}_\Gamma(\gamma_1\otimes\dots\otimes\gamma_n)\Bigr)
\end{equation}
Here $\mathrm{Alt}$ is the sum over all permutations of
$\gamma_1,\dots,\gamma_n$ with signs, such that when two polyvector fields
$\gamma_i$ and $\gamma_j$ are permuted, the sign
$(-1)^{(\deg\gamma_i+1)\cdot(\deg\gamma_j+1)}$ appears;
$G_{n,2n-3}$ stands for the set of the connected admissible graphs with $n$
vertices and $2n-3$ edges.

\begin{example}
$\mathcal{L}_2(\gamma_1\wedge\gamma_2)=\{\gamma_1,\gamma_2\}$ is the
Schouten-Nijenhuis bracket.
\end{example}

\begin{remark}
In what follows, we consider dg Lie algebras, or, more generally,
$L_\infty$ algebras. There is an ambiguity with the sign conventions. Indeed, for a
dg Lie algebra there are two possible ways to define the skew-commutativity of
the bracket: the first is $[a,b]_1=(-1)^{\deg a\deg b +1}[b,a]_1$,
and the second is $[a,b]_2=(-1)^{(\deg a+1)(\deg b+1)}[b,a]_2$. In
what follows we stick to the second definition of skew-commutativity.
\end{remark}

Recall that an $L_\infty$ algebra structure on a $\mathbb{Z}$-graded
vector space $\g$ is a coderivation $Q$ of degree +1 of the free
coalgebra $S(\g[1])$ such that $Q^2=0$. If $\g$ is a Lie algebra,
such $Q$ is given by the differential in the Chevalley-Eilenberg chain complex of the Lie
algebra $\g$.

In coordinates, an $L_\infty$ structure on $\g$ is given by a collection of
maps (the "Taylor components" of the $L_\infty$ structure)
\begin{equation}\label{k2.2.5}
\mathcal{L}_n\colon \Lambda^n\g\to\g[2-n],\, n\ge 1
\end{equation}
which satisfy, for each $N\ge 1$, the following quadratic equation:
\begin{equation}\label{k2.2.6}
\mathrm{Alt}_{g_1,\dots,g_N}\sum_{a+b=N+1,\,a,b\ge
1}\pm\frac1{a!b!}\mathcal{L}_b\Bigl((\mathcal{L}_a(g_1\wedge\dots\wedge
g_a)\wedge g_{a+1}\wedge\dots\wedge g_{N})\Bigr)=0
\end{equation}

Now turn back to the polyvector field
$\mathcal{L}_\Gamma(\gamma_1\otimes\dots\otimes\gamma_n)$, for an admissible graph $\Gamma$. This
polyvector field has degree $\sum_{i=1}^n\deg\gamma_i+n-\sharp
E(\Gamma)-1$. When $\sharp E(\Gamma)=2n-3$, this degree is
$\sum\deg\gamma_i-n+2$, what agrees with the shift of degrees in (\ref{k2.2.5}).
\begin{theorem}\label{thexotic}
The maps $\mathcal{L}_n\colon\Lambda^nT_\fin(V)\to T_\fin(V)[2-n]$
are the Taylor components of an $L_\infty$ algebra structure on
$T_\fin(V)$.
\begin{proof}
Consider the relation (\ref{k2.2.6}) for some fixed $N$. It can be
rewritten as
\begin{equation}\label{k2.2.7}
\text{the l.h.s. of }(\ref{k2.2.6})=\sum_{\Gamma\in
G_{N,2N-4}}c_\Gamma\cdot\mathrm{Alt}_{\gamma_1,\dots,\gamma_N}\mathcal{L}_\Gamma(\gamma_1\otimes\dots\otimes\gamma_N)=0
\end{equation}
where the summation is taken over all the connected admissible graphs
with $N$ vertices and $2N-4$ edges (for 1 edge less than in
(\ref{k2.2.4})), and $c_\Gamma$ are some (real) numbers. We need to
prove that $c_\Gamma=0$ for each $\Gamma\in G_{N,2N-4}$.

For, consider the integral
\begin{equation}\label{k2.2.8}
\int_{C_{N,\Gamma}}d\bigl(\bigwedge_{e\in E(\Gamma)}\phi_e\bigr)
\end{equation}
This integral is clearly equal to 0, because all 1-forms $\phi_e$
are closed (moreover, they are exact). Now we want to apply the
Stokes' formula. For this we need to construct the compactifications
$\overline{C}_{n,\Gamma}$ of the spaces $C_{n,\Gamma}$ which is a
smooth manifold with corners, and such that the forms $\phi_e$ can
be extended to a smooth forms on $\overline{C}_{n,\Gamma}$. It can
be done in the standard way, see [K97, Sect. 5]. Here we describe
the boundary strata of codimension 1 which are the only strata which
contribute to the integrals we consider. Here is the list of the
boundary strata of codimension 1:
\begin{itemize}
\item[T1)] some $S$ points $p_{i_1},\dots,p_{i_S}$ among the $n$ points approach each other,
such that $2\le \sharp S\le n-1$; in this case let $\Gamma_1$ be the
restriction of the graph $\Gamma$ into these $S$ points, and let
$\Gamma_2$ be the graph obtained from contracting of the $S$
vertices into a single new vertex. Thus, $\Gamma_1$ has $S$
vertices, and $\Gamma_2$ has $n-S+1$ vertices. In this case the
boundary stratum is isomorphic to $C_{S,\Gamma_1}\times
C_{n-S+1,\Gamma_2}$;
\item[T2)] a point $q$ connected by an edge $\overrightarrow{pq}$ with a point $p$
approaches the horizontal line passing through the point $p$.
\end{itemize}

We continue:
\begin{equation}\label{k2.2.9}
0=\int_{C_{N,\Gamma}}d\bigl(\bigwedge_{e\in
E(\Gamma)}\phi_e\bigr)=\int_{\overline{C}_{N,\Gamma}}d\bigl(\bigwedge_{e\in
E(\Gamma)}\phi_e\bigr)=\int_{\partial
\overline{C}_{N,\Gamma}}\bigwedge_{e\in E(\Gamma)}\phi_e
\end{equation}
Only the boundary strata of codimension 1 do contribute to the
r.h.s. integral. The strata of type T2) do not contribute because
the form $\bigwedge_e\phi_e$ vanishes there. We can therefore
consider only the strata of type T1).

For these strata we have the following {\it factorization
principle:} it says that the integral over a stratum $T$ of type T1)
is the product:
\begin{equation}\label{k2.2.10}
\int_T\bigwedge_{e\in\Gamma}\phi_e=\Bigl(\int_{C_{S,\Gamma_1}}\bigwedge_{e\in
\Gamma_1}\phi_e\Bigr)\times\Bigl(\int_{C_{n-S+1,\Gamma_2}}\bigwedge_{e\in
\Gamma_2}\phi_e\Bigr)
\end{equation}
The same factorization holds therefore for the weights $W_\Gamma$.

We get the following identity:
\begin{equation}\label{k2.2.11}
0=\int_{\partial T1\overline{C}_{N,\Gamma}}\bigwedge_{e\in
E(\Gamma)}\phi_e=\sum_{T\in \partial
T1}\Bigl(\int_{C_{S,\Gamma_1}}\bigwedge_{e\in
\Gamma_1}\phi_e\Bigr)\times\Bigl(\int_{C_{n-S+1,\Gamma_2}}\bigwedge_{e\in
\Gamma_2}\phi_e\Bigr)
\end{equation}
where the strata $T$ come with its orientation.

The summands in the r.h.s. are in 1-1 correspondence with the
summands in (\ref{k2.2.6}) which contribute to $c_\Gamma$.
Therefore, all $c_\Gamma=0$.
\end{proof}
\end{theorem}

It follows from Lemma \ref{lemma1721} that the $L_\infty$ structure given by Theorem \ref{thexotic} has only components of 
even $\ge 2$ degrees:
$\mathcal{L}_2,\mathcal{L}_4,\mathcal{L}_6,\dots$. On the other
hand, Example \ref{example1722} shows that the higher components
$\mathcal{L}_{2n}$, $n\ge 2$, are nonzero.

We denote the $L_\infty$ algebra $T_\fin(V)$ with the constructed
$L_\infty$ structure by $T_\fin^{\mathcal{L}}(V)$, and call it the {\it exotic} $L_\infty$ structure on $T_\fin(V)$.

\subsection{Infinite-dimensional formality}
Now we are ready to state the main result of this paper:
\begin{MainTheor}\label{theoremmain}
Let $V=\bigoplus_{i\in\mathbb{Z}_{\ge 0}}V_i$ be a non-negatively
graded vector space with $\dim V_i<\infty$. Then there is an
$L_\infty$ quasi-isomorphism from the exotic $L_\infty$ algebra
$T_\fin^{\mathcal{L}}(V)$ constructed in Section 2.2 to the
Hochschild complex $\Hoch^\mb_\fin(S(V^*))$ with the Gerstenhaber
bracket. Its first Taylor component is given by the
Hochschild-Kostant-Rosenberg map.
\end{MainTheor}
The proof of a more precise statement is given in Section 3, see Theorem \ref{mtheorembis}. Here we discuss some direct consequences.

There is an immediate corollary, obtained by comparison of  Theorem \ref{theoremmain} with the Kontsevich formality [K97]:
\begin{coroll}
Let $V$ be finite-dimensional. Then the exotic $L_\infty$ structure
$T_\fin^{\mathcal{L}}(V)$ on polyvector fields is $L_\infty$
quasi-isomorphic to the classical graded Lie algebra structure on it (given by
the Schouten-Nijenhuis bracket).
\begin{proof}
By the loc.cit., there is an $L_\infty$
quasi-isomorphism $\U\colon T_\poly(V)\to\Hoch^\mb(S(V^*))$. By our
Main Theorem, there is an $L_\infty$ quasi-isomorphism
$\mathcal{F}\colon T_\fin^{\mathcal{L}}(V)\to \Hoch^\mb(S(V^*))$. It
implies that the two $L_\infty$ structures on polyvector fields are
$L_\infty$ quasi-isomorphic.
\end{proof}
\end{coroll}
However, for an infinite-dimensional space $V$, the two structures
may be completely different (and, in fact, they are, as follows from the failure of Kontsevich formality in this case). For example, suppose we have a bivector
$\alpha\in T_\fin(V)$ for an infinite-dimensional $V$. The right
concept what is that $\alpha$ is Poisson is given by our Main
Theorem \ref{theoremmain} as follows:

\begin{equation}\label{k2.3.1}
\frac12\mathcal{L}_2(\alpha\wedge\alpha)+\frac1{24}\mathcal{L}_4(\alpha\wedge\alpha\wedge\alpha\wedge\alpha)+\dots=0
\end{equation}
For a fixed $\alpha$ which is polynomial in coordinates the sum is
actually finite. We see, in particular, that if $\alpha$ satisfies
(\ref{k2.3.1}), the bivector field $\lambda\cdot\alpha$,
$\lambda\in\mathbb{C}$, may not satisfy. That is, the equation
(\ref{k2.3.1}) is not homogeneous. We call (\ref{k2.3.1}) {\it the
quasi-Poisson equation}.

In the case when $\alpha$ is a linear bivector field the
quasi-Poisson equation (\ref{k2.3.1}) coincides with the classical
one:
\begin{lemma}
Let $\alpha$ be a linear bivector on $V$. Then the higher components
$\mathcal{L}_{2n}(\alpha^{\wedge 2n})=0$, $n\ge 2$. That is,
(\ref{k2.3.1}) is equivalent to the Poisson equation
$\{\alpha,\alpha\}=0$.
\begin{proof}
Any admissible graph with $k$ vertices which contributes to the
$L_\infty$ algebra $\mathcal{L}$ has $2k-3$ edges. If $k>2$ it
implies that there is at least one vertex with at least two incoming
edges. Then the corresponding operator
$\mathcal{L}_\Gamma(\alpha^{\wedge k})$ is zero because $\alpha$ has
linear coefficients.
\end{proof}
\end{lemma}

We expect that there are Poisson bivectors which are of degree $\ge
2$ on an infinite-dimensional vector space which {\it are impossible
to quantize} in the sense of deformation quantization. As the
condition (\ref{k2.3.1}) is non-homogeneous, our formality theorem
gives a deformation quantization only if all
$\mathcal{L}_{2n}(\alpha^{\wedge 2n})=0$, $n\ge 2$, separately, so we have a
sequence of {\it homogeneous} equations. One can say that this series
of equations gives the {\it higher obstructions} for deformation
quantization problem in the infinite-dimensional case.

\smallskip
\smallskip

Compute now the first obstruction $\mathcal{L}_4(\alpha^{\wedge
4})$. The list of all  6 admissible graphs (up to the labelling) with
4 vertices and $2\cdot 4-3=5$ edges is shown in Figure \ref{figure4}.

\sevafigc{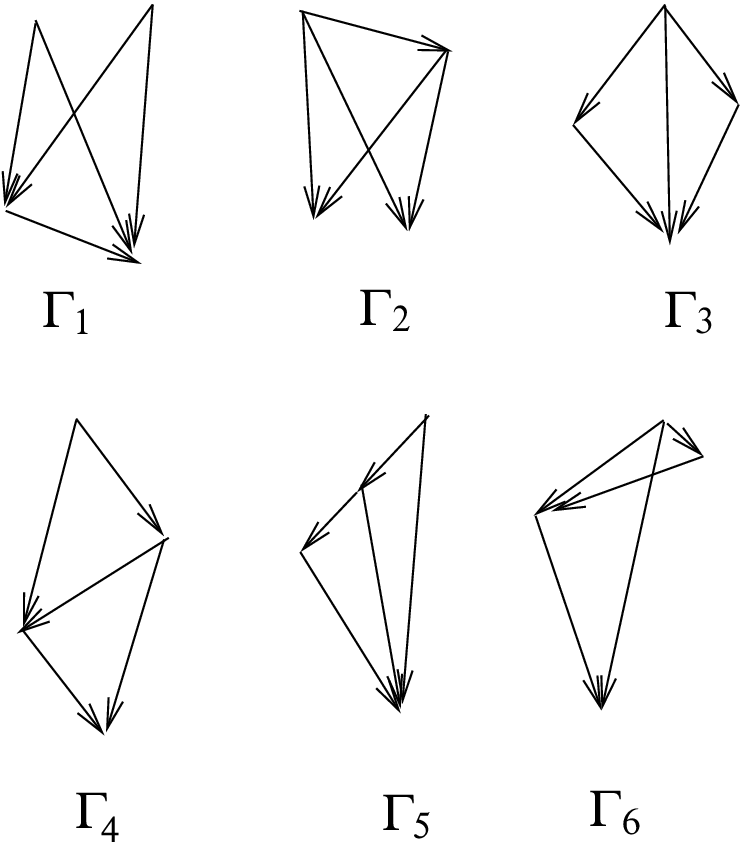}{70mm}{0}{The admissible graphs in $G_{4,5}$\label{figure4}}

Their weights are $2\times\frac{13}{12}$, $2\times\frac{13}{12}$,
$2\times\frac13$, $\frac{13}{12}$, $\frac7{12}$, $\frac7{12}$,
respectively. (Here $2\times\dots$ counts the two possible
labelings). Clearly only $\Gamma_1$, $\Gamma_4$, $\Gamma_5$
contribute to $\mathcal{L}_4(\alpha^{\wedge 4})$, because other
graphs has a vertex with 3 outgoing edges. We get the following
equation for the vanishing of the first obstruction in the
infinite-dimensional case:
$$\boxed{
\Bigl(\frac{13}{6}\mathcal{L}_{\Gamma_1}+\frac{13}{12}\mathcal{L}_{\Gamma_4}+
\frac7{12}\mathcal{L}_{\Gamma_5}\Bigr)(\alpha\wedge\alpha\wedge\alpha\wedge\alpha)=0
}$$

When $\alpha$ is a {\it quadratic} Poisson bivector, only
$\mathcal{L}_{\Gamma_4}$ survives because other graphs have a vertex
with 3 incoming edges. Then in the quadratic case the first
obstruction reads
\begin{equation}\label{quadr}
\mathcal{L}_{\Gamma_4}(\alpha\wedge\alpha\wedge\alpha\wedge\alpha)=0
\end{equation}

\section{A proof of the Main Theorem \ref{theoremmain}}
\subsection{The new propagator}\label{sectionnewprop}
Recall the configuration space $C_{n,m}$ introduced in [K97]. First of all, $\Conf_{n,m}$ is the configuration
space of $m+n$ pairwise distinct points among which $n$ belong to the open upper
half-plane $\{z\in\mathbb{C}, \mathrm{Im}z>0\}$, and the remaining $m$ belong to the real line, which
is thought on as the boundary of the open upper-half plane. Then
\begin{equation}\label{2.1.1}
C_{n,m}=\Conf_{n,m}/G_2
\end{equation}
where $G_2$ is the two-dimensional group of symmetries of the upper
half-plane of the form $\{z\mapsto az+b|, a\in\mathbb{R}_+,
b\in\mathbb{R}\}$. Here a point $z$ of the upper half-plane is
considered as a complex number with positive imaginary part.

M.Kontsevich provided loc.cit. a compactification
$\overline{C}_{n,m}$ of these spaces.  We refer
the reader to [loc.cit., Sect. 5], for a detailed description of this
compactification.  Then he constructed a top degree
differential form on $\overline{C}_{n,m}$, associated with a Kontsevich admissible graph $\Gamma$ with $n$ vertices of the first type and $m$ vertices of the second type (see [loc.cit., Sect. 6.1]). 

The top degree differential form on $\overline{C}_{n,m}$ is
constructed as follows. One firstly constructs a closed 1-form
$\phi$ on the space $\overline{C}_{2,0}$ (``the propagator"). Assume that a graph $\Gamma$ with $n+m$ vertices has exactly $2n+m-2$ oriented edges. For each edge $e$ of $\Gamma$one has the forgetful
map $t_e\colon \overline{C}_{n,m}\to \overline{C}_{2,0}$. Now if the
1-form $\phi$ is chosen, one defines the top degree form on
$\overline{C}_{n,m}$ associated with an admissible graph $\Gamma$ as
\begin{equation}\label{2.1.2}
\phi_\Gamma=\bigwedge_{e\in E(\Gamma)}t_e^*(\phi)
\end{equation}
(one should impose some order on the edges of $\Gamma$ to define the
wedge-product of 1-forms in \eqref{2.1.2}; this order is a part of data for a Kontsevich
admissible graph).

We will need a slightly refined concept of a Kontsevich admissible graph (see Section \ref{2017graphs} or [K97, Sect. 6.1]). We call this concept {\it a normalized} Kontsevich admissible graph; by the {\it normalization} we mean the condition \eqref{k3.1.1} below. 

\begin{defn}\label{normkontgraphs}{\rm
A {\it normalized Kontsevich admissible graph} is 
a Kontsevich admissible graph (see Definition \ref{kontgraphs}) such that 
the following condition holds (where for a vertex $v$ of the first type $n(v)$ denotes it labelling):
\begin{equation}\label{k3.1.1}
\text{For two vertices of the first type connected by an edge $\overrightarrow{v_1v_2}$ one has $n(v_1)<n(v_2)$}
\end{equation}
}
\end{defn}
Note that \eqref{k3.1.1} implies that {\it there are no oriented cycles between the vertices of the first type}.

\hspace{1mm}

We do not consider a choice of the orderings of the sets $Star(v)$
of outgoing edges for the vertices $v$ of $\Gamma$ as a part of the data of
an admissible graph, see Remark \ref{remark171new}. 

Recall firstly what the space $\overline{C}_{2,0}$ is. It looks like
``an eye'' (see the left picture in Figure \ref{figure5}). The two boundary lines
comes when one of the two points $z_1$ or $z_2$ approaches the real
line, which is the boundary of the upper half-plane. The circle
is the boundary stratum corresponded to the case when the points $z_1$ and $z_2$ approach each other and
are far from the real line. The role of the two points $z_1$ and
$z_2$ here is completely symmetric. Now we are going to break this
symmetry down.

Subdivide the space $\overline{C}_{2,0}$ into two subspaces, as follows.
Let the point $z_1$ be fixed. Draw the half-circle orthogonal to the
real line (a geodesic in the Poincar\'{e} model of hyperbolic
geometry) such that $z_1$ is the top point of the half-circle, see Figure \ref{figure6}. The
half-circle is a geodesic, and the group $G_2$ in (\ref{2.1.1}) is the
group of symmetries for the Poincar'{e} model of hyperbolic geometry. It proves that a
half-circle orthogonal to the real line is transformed to an analogous half-circle,  by any $g\in
G_2$. Therefore, the image of the half-circle is well-defined on the
``eye '' $\overline{C}_{2,0}$.

We draw this image in Figure \ref{figure5}, as the border line between
the light and the dark parts. We want the 1-form
$\phi$ to vanish when the oriented pair $(z_1,z_2)$ is in the light area of
Figure \ref{figure6} ($z_2$ is outside the half-circle). We contract all light
area in the left-hand side picture in Figure \ref{figure5} to a point, and get the
right-hand side picture therefrom. Here both vertices of the eye from the left-hand side
picture, and the upper boundary component of it, as
well as the upper half of the circle, are contracted to a one point. The external boundary
in the right-hand side picture is formed by the low boundary component in
the left-hand side picture. We call the space drawn in the right-hand side picture in Figure \ref{figure5} {\it the modified} (contracted) $\overline{C}_{2,0}$, and denote it by $\overline{C}_{2,0}^m$.

\begin{defn}
A {\it modified angle function} is any map from the modified (contracted)
$\overline{C}_{2,0}$ to the circle unit $S^1$ such that the internal
circle is mapped isomorphically to
$S^1$, as the Euclidean ${\mathrm{angle}}$, and the external boundary
component is also mapped to $S^1$ (in the homotopically unique way, with
the period $\pi$).
\end{defn}
Here the angle in the internal circle is the Euclidian angle when
$z_2$ approaches $z_1$. The period of this angle inside the
half-circle is $\pi$, not $2\pi$.

There is a modified Kontsevich's harmonic angle function, which provides an example of such a map, see
(\ref{new1}) below.

Let us recall that {\it the Kontsevich angle function}, introduced in [K97, Sect. 6.2], is a
map of the eye $\overline{C}_{2,0}$ (drawn in the left picture of Figure \ref{figure5}) to $S^1$ such that the inner
circle is mapped isomorphically as the angle function, the upper boundary is
contracted to a point, and the lower boundary is mapped with period
$2\pi$.

Let us stress a difference between the two definitions: in our
case, any angle function $\theta$ is a function in a proper sense,
whence in Kontsevich's definition it is a multi-valued function
defined only up to $2\pi$. Therefore, in our case, the de Rham
derivative $\phi=d\theta$ is an exact 1-form, whence in 
Kontsevich's case it is only closed. 
What also makes this difference essential
is the observation
that in our case the propagator $\phi$ is not a smooth function on
the manifold with corners $\overline{C}_{2,0}$. It will be the
main source of problems for the proof of the $L_\infty$ identities using Stokes' formula in the next Subsection.

Define the {\it weight of an admissible graph} $\Gamma$ as
\begin{equation}\label{2.1.7}
W_\Gamma=\frac1{\pi^{2n+m-2}}\int_{\overline{C}_{n,m}^+}\bigwedge_{e\in
E(\Gamma)}f_e^*\phi
\end{equation}
where  $\phi$ is a modified angle function, and $f_e$ is the forgetful map $f_e\colon \overline{C}_{n,m}\to \overline{C}_{2,0}$ associated with an edge $e\in
E(\Gamma)$.

This definition although being correct has a small drawback, due to the non-smoothness of the
1-form $\phi$ on $\overline{C}_{2,0}$. Potentially it may make troubles in proofs using the Stokes formula on manifold with corners. We encourage to the reader not be confused by this drawback now, and to wait untill Section 3.2.2, when we provide a better reformulation of the
definition of weights, in (\ref{2.2.5}).

Note that, although the 1-form $\phi$ is exact, the integrals
$W_\Gamma$ are in general nonzero, because the spaces
$\overline{C}_{n,m}$ are manifolds (with corners) with boundary.

\sevafigc{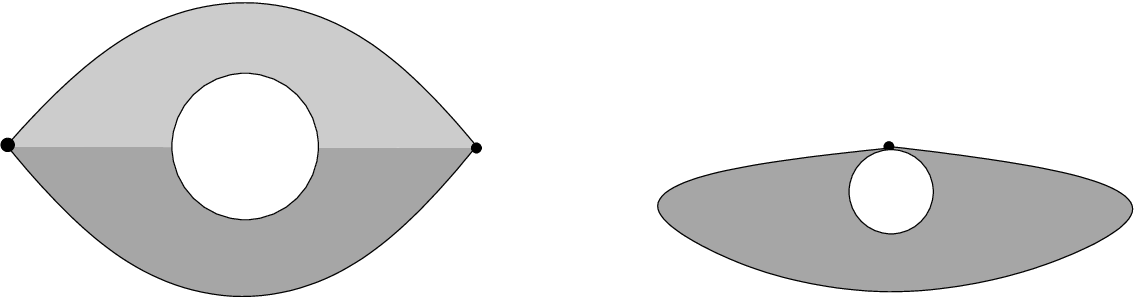}{110mm}{0}{The Kontsevich's configuration space
$\overline{C}_{2,0}$ (left) and our modified space $\overline{C}_{2,0}^m$ (right)\label{figure5}}

\sevafigc{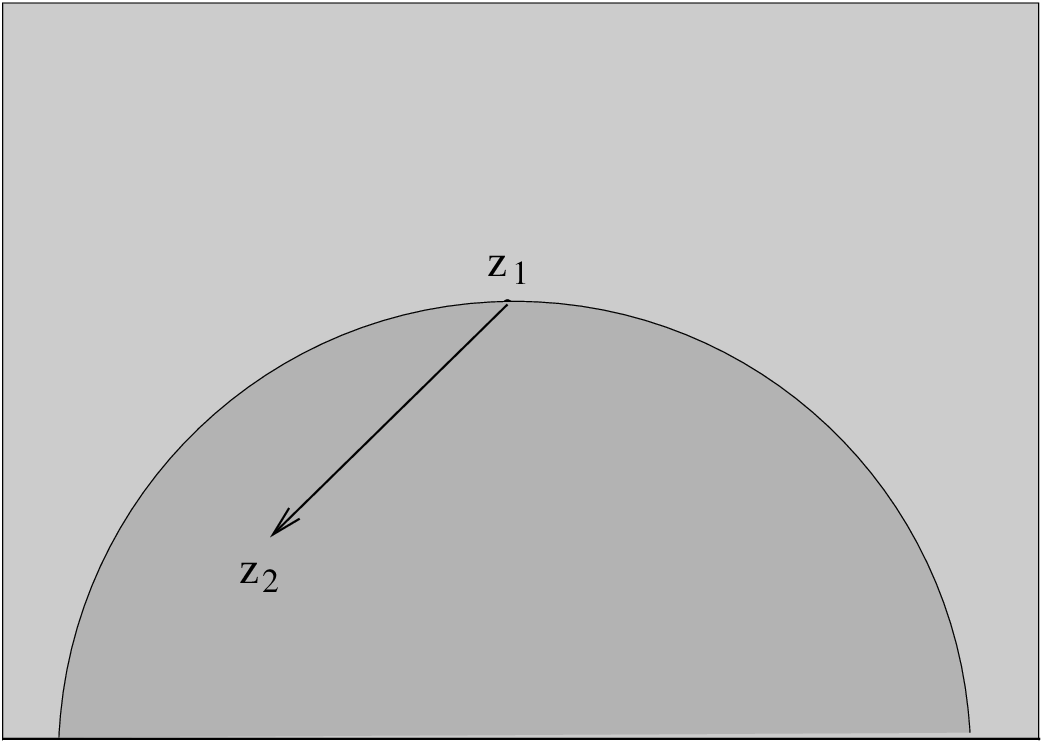}{80mm}{0}{The height ordering of two points,
$z_2\le z_1$\label{figure6}}

\subsection{The proof}
\subsubsection{Statement of the result}
Let $V$ be a $\mathbb{Z}_{\ge 0}$-graded vector space over
$\mathbb{C}$ with finite-dimensional components $V_i$.

Recall our concept of normalized Kontsevich admissible graphs, see Definition \ref{normkontgraphs}.
 Let $G_{n,m,2n+m-2}$ be the set of all connected
admissible graphs with $n$ vertices of the first type, $m$ vertices
of the second type, and $2n+m-2$ edges. Recall the polydifferential operators
$\U_\Gamma(\gamma_1\wedge\dots\wedge\gamma_n)$, $\Gamma\in G_{n,m}$,
associated with a (general) Kontsevich admissible graph $\Gamma$
and with polyvector fields $\gamma_1,\dots,\gamma_n$, see [K97, Sect.
6.3]. Define
\begin{equation}\label{2.2.1}
\F_n(\gamma_1\wedge\dots\wedge\gamma_n)=\sum_{\Gamma\in
G_{n,m}}W_\Gamma\times\U_\Gamma(\gamma_1\wedge\dots\wedge\gamma_n)
\end{equation}
where the weight $W_\Gamma$ is defined in (\ref{2.1.7}) via the
modified angle function. Consider the following two cases: either
$\Gamma$ contains an oriented cycle (and in this case clearly our
$W_\Gamma=0$), or it does not contain any oriented cycle (and in the latter case $\U_\Gamma$
is well-defined by Lemma \ref{lemma1714}). Therefore, the cochain
$\F_n(\gamma_1\wedge\dots\wedge\gamma_n)$ is well-defined.

We prove
\begin{theorem}\label{mtheorembis}
Let $V$ be a $\mathbb{Z}_{\ge 0}$-graded vector space over
$\mathbb{C}$ with finite-dimensional graded components $V_i$. Then
the maps $\F_n$ are well-defined and are the Taylor components of
an $L_\infty$ quasi-isomorphism $\F\colon T_\fin^\mathcal{L}(V)\to
\Hoch^\mb_\fin(S(V^*))$. (Here in the l.h.s. $T_\fin(V)^\mathcal{L}$ is the exotic $L_\infty$ algebra 
introduced in Section 2). The first Taylor component $\F_1$ of the $L_\infty$ quasi-isomorphism $\F$ is the
Hochschild-Kostant-Rosenber map (\ref{1.3.1}).
\end{theorem}
Theorem \ref{mtheorembis} provides an explicit form of the $L_\infty$ quasi-isomorphism which existence is stated in Main Theorem \ref{theoremmain}, and, therefore, it implies the latter Theorem. Below we prove Theorem \ref{mtheorembis}.

\subsubsection{Configuration spaces}
To prove Theorem \ref{mtheorembis}, we would like to apply the well-known argument with Stokes' formula for manifolds with corners, see e.g. [K97, Sect. 6.4].
The main trouble, which makes it impossible to apply this argument straightforwardly, is that the de Rham
derivative of the modified angle function is {\it not} a smooth
differential form on the compactified space $\overline{C}_{2,0}$. There could be (at least) two different ways to overcome this problem. The first one
is to subdivide the configuration spaces $\overline{C}_{n,m}$ such
that the wedge-product of the angle 1-forms would be a smooth
differential form on the subdivision. In particular, it is clear how
to subdivide $\overline{C}_{2,0}$: the subdivision is shown in
Figure \ref{figure7}. This approach, however, wouldn't work, because the
subdivided spaces are not manifold with corners anymore, and one can
not apply the Stokes formula for them.
\sevafigc{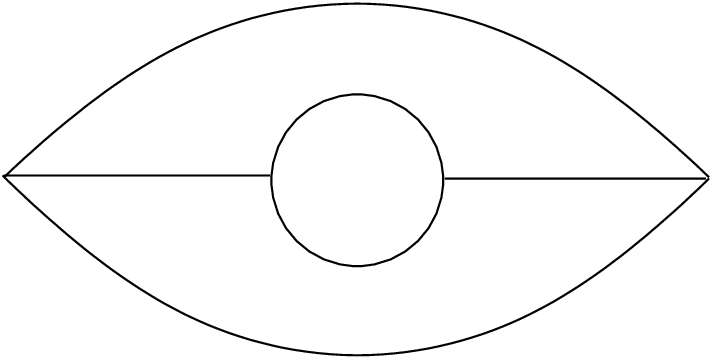}{70mm}{0}{A possible subdivision
$\overline{C}_{2,0}$ (which is not a manifold with corners)\label{figure7}}

The second way is to define another configuration spaces, cutting
the area $z_1\le z_2$ off if there is an edge
from $z_1$ to $z_2$ (see Figure \ref{figure6}) . This solution would complicate the constructions,  as 
{\it the configuration spaces would depend on a graph $\Gamma$}, that is, for each
$\Gamma$ we would have its own configuration space. However, we will show in the rest of the paper that this solution works.

Let $\Gamma$ be an normalized Kontsevich admissible graph with $n$ vertices of the first
type and $m$ vertices of the second type. Recall that
all vertices of $\Gamma$ of the first type are 
labeled as $\{1,2,\dots,n\}$, and all vertices of the second type
are ordered and labeled as
$\{\overline{1},\overline{2},\dots,\overline{m}\}$. Recall also that \eqref{k3.1.1} is assumed. Define the
configuration space $C_{n,m,\Gamma}$ as follows:
\begin{equation}\label{2.2.2}
\begin{aligned}
\ &C_{n,m,\Gamma}=\{z_1,\dots,z_n\in\mathcal{H};\, t_1,\dots, t_n\in
\mathbb{R}|\\
&\text{$z_i\ne z_j$ for $i\ne j$, $t_1<\dots<t_m$,  $z_j\le z_i$ if
$\overrightarrow{(i,j)}$ is an edge in $\Gamma$}\\
&\text{and $t_j\le z_i$ if $\overrightarrow{(i,\overline{j})}$ is an
edge in $\Gamma$}\}/G_2
\end{aligned}
\end{equation}
Here $\mathcal{H}=\{z\in \mathbb{C}, \mathrm{Im}z>0\}$ is the upper
half-plane, and the relations $z_j\le z_i$ and $t_j\le z_i$ in the
r.h.s. are understood in the sense of the ordering with half-circle,
shown in Figure \ref{figure6}. The group $G_2$ in the r.h.s. is $G_2=\{z\mapsto az+b, a\in \mathbb{R}_{>0},
b\in\mathbb{R}\}$.

Note that the graph $\Gamma$ in the definition of
$C_{n,m,\Gamma}$ may have an arbitrary (not necessarily $2n+m-2$)
number of edges. If $\Gamma$ has no edges at all, we recognize the
Kontsevich's original space $C_{n,m}$.

It is easy to construct a Kontsevich-type compactification
$\overline{C}_{n,m,\Gamma}$ which is a manifold with corners, with
projections $p_{\Gamma,\Gamma^\prime}\colon
\overline{C}_{n,m,\Gamma}\to \overline{C}_{n^\prime,
m^\prime,\Gamma^\prime}$ defined when $\Gamma^\prime$ is a subgraph
of $\Gamma$. The space $\overline{C}_{2,0,\Gamma_0}$, where
$\Gamma_0$ is just one oriented edge connecting point 1 with point
2, is shown in Figure \ref{figure8}.

\sevafigc{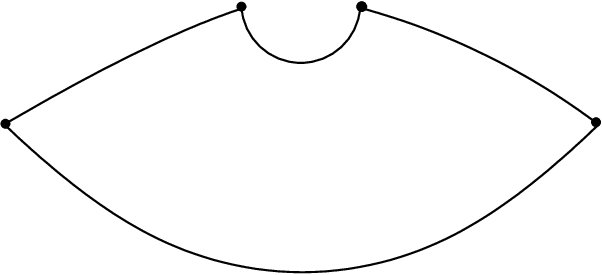}{70mm}{0}{The "Eye"
$\overline{C}_{2,0,\Gamma_0}$\label{figure8}}

Here the lower component of the boundary comes when the point $z_2$
(see Figure \ref{figure6}) approaches the real line, the left (corresp., the right)
upper boundary component is corresponded to the configurations when the point $z_2$ approaches the
left (corresp., the right) part of the geodesic half-circle in Figure \ref{figure6},
and the third upper boundary component (the half-circle) comes when
$z_2$ approaches $z_1$ inside the geodesic half-circle.

We can upgrade the definition of the modified angle function from this
new point of view:

\begin{defn}
A {\it modified angle function} is a continuous map $\theta\colon
\overline{C}_{2,0,\Gamma_0}\to S^1$ such that
$\theta$ is given by the Euclidean angle varying from $-\pi$ to 0 on the
upper half-circle, and $\theta$ contracts the two other upper
boundary components to a point $0\in S^1$.
\end{defn}

{\it An example} of the modified angle function is the doubled
Kontsevich's harmonic angle:

\begin{equation}\label{new1}
f(z_1,z_2)=\frac1{i}\mathrm{Log}\frac{(z_1-z_2)(z_1-\overline{z}_2)}{(\overline{z}_1-z_2)(\overline{z}_1-\overline{z}_2)}\text{
where $z_1\ge z_2$ in the sense of Figure \ref{figure6}}
\end{equation}
It is clear that the function $f(z_1,z_2)$ defined in this way is
well-defined when $z_1\ge z_2$ and is equal to 0 on the "border"
circle, see Figure \ref{figure6}.

Now we define for an edge $e$ of an admissible graph $\Gamma$ the
1-form $\phi_e$ on $\overline{C}_{n,m,\Gamma}$ as
$p_{\Gamma,\Gamma_0}^*(d\theta)$, where $\Gamma_0$ is the graph with
two vertices and one edge $e$. Finally, we give a rigorous
definition of the weight:

\begin{equation}\label{2.2.5}
W_\Gamma=\frac1{\pi^{2n+m-2}}\int_{\overline{C}_{n,m,\Gamma}}\bigwedge_{e\in
E(\Gamma)}\phi_e
\end{equation}

Let us describe the boundary strata of codimension 1 of
$\overline{C}_{n,m,\Gamma}$. The boundary strata will be expressed in terms of spaces $\overline{C}_{n^\prime,m^\prime,\Gamma^\prime}$, as well as of spaces
$C_{n,\Gamma}$, introduced in Section 2.

Below is a complete list of the types of boundary strata of codimension 1 in
$\overline{C}_{n,m,\Gamma}$:
\begin{itemize}
\item[S1)] some points $p_1,\dots,p_S\in\mathcal{H}$, $S\ge 2$, approach each other
and remain far from the geodesic half-circles of any point $q\ne
p_1,\dots, p_S$ such that there is an edge $\overrightarrow{qp_i}$, for some $1\le i\le S$. In this case we get the boundary stratum of codimension 1
isomorphic to $C_{n-S+1,m,\Gamma_1}\times C_{S,\Gamma_2}$ where
$\Gamma_2$ is the subgraph of $\Gamma$ of the edges connecting the
point $p_1,\dots, p_S$ with each other, and $\Gamma_1$ is obtained
from $\Gamma$ by collapsing the graph $\Gamma_2$ into a (new) single vertex;
\item[S2)] some points $p_1,\dots,p_S\in\mathcal{H}$ and some points
$q_1,\dots, q_R\in\mathbb{R}$, $2S+R\ge 2$, $2S+R\le 2m+n-1$,
approach each other and a point of the real line, but are far from
the geodesic half-circle of any other point connected with these
points by an edge. In this case we get a boundary stratum of
codimension 1 isomorphic to $C_{S,R,\Gamma_1}\times
C_{n-S,m-R+1,\Gamma_2}$ where $\Gamma_1$ is the graph formed from
the edges between the approaching points, and $\Gamma_2$ is obtained
from $\Gamma$ by collapsing the subgraph $\Gamma_1$ into a new
vertex of the second type;
\item[S3)] some point $p$ is placed on  the geodesic
half-circle of exactly one point $q\ne p$ which is far from $p$,
such there is an edge $\overrightarrow{qp}$. These boundary
strata will be irrelevant for us because they do not contribute to
the integrals we consider.
\end{itemize}

We are ready to prove Theorem \ref{mtheorembis}.

\subsubsection{Application of the Stokes' formula and the boundary
strata} We need to prove the $L_\infty$ morphism quadratic relations
on the maps $\F_n$. Recall, that for each $k\ge1$ and for any
polyvector fields $\gamma_1,\dots,\gamma_k\in
T_\fin^{\mathcal{L}}(V)$ it is the relation:
\begin{equation}\label{2.2.10}
\begin{aligned}
\ &d_\Hoch(\F_k(\gamma_1\wedge\dots\wedge\gamma_k))+\\
&\sum_{2\le N\le k}\frac1{k!(N-k)!}\sum_{\sigma\in S_k}\pm\F_{k-N}\bigl(\mathcal{L}_N(\gamma_{\sigma(1)}\wedge\dots\wedge\gamma_{\sigma(N)})\wedge\gamma_{\sigma(k-N+1)}\wedge\dots\wedge\gamma_{\sigma(k)}\bigr)+\\
&\frac12\sum_{a,b\ge 1,a+b=k}\frac1{a!b!}\sum_{\sigma\in
S_k}\pm[\F_a(\gamma_{\sigma(1)}\wedge\dots\wedge\gamma_{\sigma(a)}),\F_b(\gamma_{\sigma(a+1)}\wedge\dots\wedge\gamma_{\sigma(k)})]=0
\end{aligned}
\end{equation}
The l.h.s. of (\ref{2.2.10}) is a sum over the admissible graphs
$\Gamma^\prime\in G_{n,m,2n+m-3}$ with $n$ vertices of the first
type, $m$ vertices of the second type, and $2n+m-3$ egdes (that is,
having for 1 edge less than the graphs contributing to $\F_n$). This sum
is of the form:
\begin{equation}\label{2.2.11}
\text{l.h.s}=\sum_{\Gamma^\prime\in
G_{n,m,2n+m-3}}\alpha_{\Gamma^\prime}\U_{\Gamma^\prime}
\end{equation}
where $\alpha_{\Gamma^\prime}$ are some complex numbers, and
$\U_{\Gamma^\prime}$ are the Kontsevich's polydifferential operators
from [K97]. {\it It is clear that all $\U_{\Gamma^\prime}$ do not
contain any oriented cycle}.

Our goal is to prove that all numbers $\alpha_{\Gamma^\prime}=0$.

Each $\alpha_{\Gamma^\prime}$ is a quadratic-linear combination of
our weights $W_{\Gamma}$ for $\Gamma\in G_{n,m,2n+m-2}$. 

There is a
general construction producing identities on quadratic-linear combinations of $W_\Gamma$, as follows:

Consider some $\Gamma^\prime\in G_{n,m,2n+m-3}$. With $\Gamma$ is associated
a differential form $\bigwedge_{e\in E(\Gamma^\prime)}\phi_e$ on
$C_{n,m,\Gamma^\prime}$, as above. Consider
\begin{equation}\label{2.2.12}
\int_{C_{n,m,\Gamma^\prime}}d\bigl(\bigwedge_{e\in
E(\Gamma^\prime)}\phi_e\bigr)
\end{equation}
This expression is 0 because the form $\bigwedge_{e\in
E(\Gamma^\prime)}\phi_e$ is closed (it is, moreover, exact). Next, by
the Stokes theorem, we have
\begin{equation}\label{2.2.13}
0=\int_{C_{n,m,\Gamma^\prime}}d\bigl(\bigwedge_{e\in
E(\Gamma^\prime)}\phi_e\bigr)=\int_{\partial
C_{n,m,\Gamma^\prime}}\bigl(\bigwedge_{e\in
E(\Gamma^\prime)}\phi_e\bigr)
\end{equation}
The only strata of codimension 1 in
$\partial{C_{n,m,\Gamma^\prime}}$ do contribute to the r.h.s. They
are given by the list S1)-S3) in Section 3.2.2.

The strata of type S3) clearly do not contribute because of our boundary
conditions on the modified angle function.

For the strata of types S1) and S2) one has the following {\it factorization
property}: the integral over the product is equal to the product of
integrals. It follows from the factorization of differential forms on the product of corresponding spaces.

It makes possible to prove the following key-lemma:

\begin{klemma}\label{klemma17}
For any normalized Kontsevich admissible graph $\Gamma^\prime\in G_{n,m,2n+m-3}$, the coefficient
$\alpha_{\Gamma^\prime}$ in (\ref{2.2.11}) is equal to
$\int_{\partial C_{n,m,\Gamma^\prime}}\bigl(\bigwedge_{e\in
E(\Gamma^\prime)}\phi_e\bigr)$ (which is zero by the Stokes'
formula).
\end{klemma}
Theorem \ref{mtheorembis} follows directly from this Key-Lemma \ref{klemma17}.

To prove Key-Lemma \ref{klemma17}, we express
\begin{equation}\label{2.2.14}
\int_{\partial C_{n,m,\Gamma^\prime}}\bigl(\bigwedge_{e\in
E(\Gamma^\prime)}\phi_e\bigr)=\int_{\partial_{S1}
C_{n,m,\Gamma^\prime}}\bigl(\bigwedge_{e\in
E(\Gamma^\prime)}\phi_e\bigr)+\int_{\partial_{S2}
C_{n,m,\Gamma^\prime}}\bigl(\bigwedge_{e\in
E(\Gamma^\prime)}\phi_e\bigr)+\int_{\partial_{S3}
C_{n,m,\Gamma^\prime}}\bigl(\bigwedge_{e\in
E(\Gamma^\prime)}\phi_e\bigr)
\end{equation}
the integral over the boundary $\partial C_{n,m,\Gamma^\prime}$ as
the sum of the integrals over the three types of boundary strata
S1)-S3) of codimension 1.

As we have already noticed, the integral $\int_{\partial_{S3}
C_{n,m,\Gamma^\prime}}\bigl(\bigwedge_{e\in
E(\Gamma^\prime)}\phi_e\bigr)=0$, due to the boundary conditions for
the propagator.

The summand $\int_{\partial_{S2}
C_{n,m,\Gamma^\prime}}\bigl(\bigwedge_{e\in
E(\Gamma^\prime)}\phi_e\bigr)$ corresponds exactly to the first and
to the third summands of the l.h.s. of (\ref{2.2.10}) containing the
Hochschild differential and the Gerstenhaber bracket, by the
factorization property.

It remains to associate a boundary stratum of type S1) with a summand of the second
summand of (\ref{2.2.10}), containing operations $\mathcal{L}_N$. It
is clear that the part of the second summand of (\ref{2.2.10}) for a
fixed $N$ is in 1-to-1 correspondence with that summands in
$\int_{\partial_{S1} C_{n,m,\Gamma^\prime}}\bigl(\bigwedge_{e\in
E(\Gamma^\prime)}\phi_e\bigr)$ where $S=N$ points in the upper
half-plane approach each other.

Theorems \ref{mtheorembis} and Theorem \ref{theoremmain} are proven.

\qed

\begin{remark}{\rm
In [K97], where computations of this type are originated from, the strata analogous to our
strata of type $S1)$, contribute zero integrals, for $N>2$. This is why
M.Kontsevich gets his ``pure'' formality theorem there. This vanishing is
proven by a quite non-trivial computation, see loc.cit., Sect. 6.6. We saw in
Examples \ref{example1721} and \ref{example1722} that the analogous vanishing fails in the framework of our
construction.}
\end{remark}

\appendix

\section{The failure of the Kontsevich formality for $\Hoch^\mb_{\fin}(S(V^*))$ for an infinite-dimensional $V$}
In this Appendix,  we provide an argument which shows that, for a general infinite-dimensional $V$, the formality of the dg Lie algebra $\Hoch^\udot_{\fin}(S(V^*))$ fails.

More precisely, we prove the following statement:
\begin{theorem}\label{theorema1}
Let $W$ be a finite-dimensional vector space over a field $\k$ of characteristic 0, and $L=\Lie(W)$ be the free Lie algebra generated by $W$. Denote by $V$ the underlying graded space of $L$, with the inherited grading:
$$
V=\oplus_{n\ge 1}V_n,\ \ V_n=\Lie_n(W)
$$
where $\Lie_n(W)$ is the subspace of $\Lie(W)$ generated by the length $n$ Lie monomials. Then there does not exist any $\gl(W)$-equivariant $L_\infty$ morphism $$\mathcal{U}\colon T_{\fin}^\udot(V^*)[1]\to\Hoch^\udot_{\fin}(S(V))[1]$$
such that the first Taylor component $\mathcal{U}_1$ is the Hochschild-Kostant-Rosenberg map.
\end{theorem}

The proof goes as follows. 

{\it Step 1.} Assume that such an $L_\infty$ morphism $\mathcal{U}$ exists. Take the linear Poisson bivector $\alpha$ in $T_\fin(V^*)[1]$ corresponded to the free Lie algebra structure on $V=\Lie(W)$ (the Kostant-Kirillov bivector). Consider $\mathcal{U}_*(\hbar\alpha)$, which is a Maurer-Cartan element in $\Hoch_\fin(S(V))[[\hbar]][1]$. Explicitly, one has:
\begin{equation}\label{appa1}
\mathcal{U}_*(\hbar\alpha)(f,g)=f\cdot g+\hbar\mathcal{U}_1(\alpha)(f\otimes g)+\frac12\hbar^2\mathcal{U}_2(\alpha,\alpha)+\dots
\end{equation}
We prove the following
\begin{prop}\label{propa1}
The star-algebra $S(V)_*[[\hbar]]$ is isomorphic to the universal enveloping algebra $\mathcal{U}(\hbar\Lie(W))[[\hbar]]$ \text{$($to the universal enveloping algebra of the Lie algebra} \text{ whose bracket is $\hbar[-,-]$, where $[-,-]$ is the Lie bracket of $\Lie(W))$}.
\end{prop}
Note that the statement analogous to this one fails for a general $L_\infty$ map $T_\poly(\g^*)[1]\to\Hoch^\udot(S(\g))[1]$ and general Lie algebra $\g$, even for the case of a finite-dimensional Lie algebra $\g$.

In particular, M.Kontsevich provided a rather specific argument in [K97, Sect. 8.3.1], to show the claim for the $L_\infty$ morphism constructed in loc.cit. 

\begin{proof}
We make use the obstruction theory introduced in [GM, Sect. 2.6]. The idea is to show that {\it any} two Maurer-Cartan elements in $T_\fin(V^*)[[\hbar]][1]$ of the form
\begin{equation}\label{appa2}
\hbar\alpha +\hbar^2\alpha_1+\hbar^3\alpha_2+\dots
\end{equation}
(where $\alpha$ is the fixed Kostant-Kirillov bivector)
are ``connected'' by a gauge transformation $\exp(v)$ corresponded to a vector field of the form
\begin{equation}
v=\hbar v_1+\hbar^2 v_2+\dots
\end{equation}

We use the well-known fact that any $L_\infty$ quasi-isomorphism of dg Lie algebras induces an isomorphism on $\pi_0(-)$ of the corresponding Deligne grouppoids. In particular, the star-product on $S(V)[[\hbar]]$, corresponded to  $\mathcal{U}(\hbar\Lie(W))[[\hbar]]$ by the PBW theorem, gives a Maurer-Cartan element of the form \eqref{appa2}, while the star-product \eqref{appa1} is corresponded to the simplest element of the form \eqref{appa2}, the one with $\alpha_1=\alpha_2=\dots=0$. It is enough to prove that any (and, in particular, these two) Maurer-Cartan elements are gauge equivalent.

In fact, the problem can be reformulated in terms of the localized dg Lie algebra $T_\loc(V^*)[[\hbar]][1]=(T_\fin(V^*)[[\hbar]][1],d=\ad\alpha)$ (in this localized dg Lie algebra, we consider the Maurer-Cartan elements of the form $\hbar\alpha_1+\hbar^2\alpha_2+\dots$, with 0 coefficient at $\hbar^0$).

One easily sees that there is an isomorphism of sets:
\begin{equation}\label{appa33}
\begin{aligned}
\ &\bigg\{\hbar\alpha+\hbar^2\alpha_1+\hbar^3\alpha_2+\dots\in \mathrm{MC}(T_\fin(V^*)[[\hbar]][1])\bigg\}/\big\{\text{  gauge actions $\exp(\hbar v_1+\hbar^2 v_2+\dots)$}\big\}\simeq\\
&\bigg\{\hbar^1\alpha_1+\hbar^2\alpha_2+\dots\in \mathrm{MC}(T_\loc(V^*)[[\hbar]][1])\bigg\}/\big\{\text{  gauge actions $\exp(\hbar v_1+\hbar^2 v_2+\dots)$}\big\}
\end{aligned}
\end{equation}
We apply the obstruction theory for extension of a gauge equivalence from [GM, Prop. 2.6(2)] for the scalars extension from $\mathbb{C}[\hbar]/\hbar^n$ to $\mathbb{C}[\hbar]/h^{n+1}$, to the dg Lie algebra $T_\loc(V^*)[[\hbar]][1]$ (corresponded to the r.h.s. of \eqref{appa33}). We need to show that the obstructions vanish for any $n$.

Accordinding to loc.cit., the obstructions for extension of the gauge transformation, belong to $H^1(\g\otimes \mathfrak{J})$, where in our case $\g=T_\loc(V^*)[[\hbar]][1]$, and $\mathfrak{J}=\hbar^n\mathbb{C}[\hbar]/\hbar^{n+1}$. Therefore, the obstructions belong to $H^2(T_\loc(V^*)[[\hbar]],-)$.

We have:
\begin{equation}
T_\loc(V^*)[[\hbar]]=C^\udot_{\mathrm{CE}}(\Lie(W), S(\Lie W)[[\hbar]])
\end{equation}
The Chevalley-Eilenberg cohomology of a free Lie algebra vanish in all degrees except for degrees 0 and 1 (for any coefficients); in particular, $H^2_{\mathrm{CE}}(\Lie(W),-)=0$. It gives the result.

\end{proof}

{\it Step 2.}
One can specialize to $\hbar=1$.
The result of Proposition \ref{propa1} gives, under the assumption contrary to the one of the statement of Theorem \ref{theorema1}, an existence of a $\gl(W)$-equivariant $L_\infty$ morphism 
\begin{equation}\label{appa3}
\mathcal{U}_\loc\colon T_\loc(V^*)[1]\to\Hoch_\fin(\mathcal{U}(\Lie(W))[1]
\end{equation}
where $T_\loc(V^*)=C_{\mathrm{CE}}(\Lie(W), S(\Lie(W)))$.
(To localize an $L_\infty$-morohism by a Maurer-Cartan equation, we use the standard well-known formula).

The cohomology of both sides of \eqref{appa3} vanish except for the degrees 0 and 1, because $\Lie(W)$ is a free Lie algebra. Recall that the universal enveloping algebra $\mathcal{U}(\Lie(W))=T(W)$ is the free associative algebra generated by $W$.

The cohomology in degree 0 is equal to $\mathbb{C}$ and is not interesting; the cohomology in degree 1 is ``very big'' and ``very interesting''. We consider the degree 1 cohomology of both sides of \eqref{appa3} as Lie algebras. 

The existence of an $L_\infty$ morphism \eqref{appa3} implies that these degree 1 cohomology are isomorphic {\it as Lie algebras}.
We show that the latter statement is false. Consequently, the statememt of Theorem \ref{theorema1} is true.

The degree 1 cohomology is easy to find. For the r.h.s., it is 
$$
\g_2=\Der_{\Assoc}(T(W))/\Inn
$$
the Lie algebra of the outer derivations of the free associative algebra $T(W)$. For the l.h.s., it is
$$
\g_1=\Der_\Pois(S(\Lie(W)))/\Inn
$$
the Lie algebra of the outer Poisson derivations of the free Poisson algebra $\Pois(W)=S(\Lie(W))$.

Assuming the contrary to the statement of Theorem \ref{theorema1}, there exists a $\gl(W)$-equivariant isomorphism 
\begin{equation}\label{appa5}
\phi\colon \g_1\to \g_2
\end{equation}

One has:
\begin{prop}\label{propa2}
There does not exist any $\gl(W)$-equivariant isomorphism $\phi$ of Lie algebras, as in \eqref{appa5}.
\end{prop}
\begin{proof}
It is rather easy to see, in fact.

The idea is to consider the Lie subalgebra $\Der(S(W))$ of $\g_1$, as follows.

Any Poisson derivation $D$ of $S(\Lie(W))$ is uniquely determined by its restriction to the generators $W\subset \Pois(W)$. In general, it is a linear map
$$
D\colon V\to S(\Lie(W))
$$
The vector space of all such derivations contains a subspace of the derivations, given by a linear map
\begin{equation}\label{appa6}
D\colon W\to S(W)
\end{equation}
These derivations form a Lie subalgebra in $\Der_\Pois(S(\Lie(W)))$, isomorphic to the Lie algebra $\Vect(W^*)=\Der_\Comm(S(W))$
of polynomial vector fields on $W^*$. 

Consider the restriction
\begin{equation}\label{appa7}
\Der_\Comm(S(W))\to \Der_\Pois(S(\Lie(W)))/\Inn\xrightarrow{\phi}\Der_\Assoc(T(W))/\Inn
\end{equation}
It is clearly an imbedding.

\begin{lemma}\label{lemmaa1}
There does not exist any $\gl(W)$-equivariant imbedding of Lie algebras $\phi_0\colon \Der_\Comm(S(W))\to \Der_\Assoc(T(W))/\Inn$.
\end{lemma}
\begin{proof}\footnote{The author is indebted to Maxim Kontsevich for the proof of Lemma \ref{lemmaa1} given below.}
We can (not canonically) imbed $\Der_\Assoc(T(W))/\Inn\hookrightarrow \Der_\Assoc(T(W))$, as $\gl(W)$-modules. Then any imbedding $\phi_0$ gives an imbedding of $\gl(W)$-modules 
\begin{equation}\label{appa10}
\hat{\phi}_0\colon \Der_\Comm(S(W))\to\Der_\Assoc(T(W))
\end{equation}
Any derivation in $\Der_\Comm(S(W))$ is a linear combination of the homogeneous maps $W\to S^N(W)$, and any derivation in $\Der_\Assoc(T(W))$ is a linear combination of the homogeneous maps $W\to T^N(W)$. The map $\hat{\phi}_0$, being $\gl(W)$-equivariant, preserves the homegeneity degree $N$. Thus, it is given by its $\gl(W)$-equivatiant components
\begin{equation}
\hat{\phi}_0^{(N)}\colon W^*\otimes S^N(W)\to W^*\otimes T^N(W)
\end{equation}
Take $N\ll \dim W$ so that the invariants can be described by the H.Weyl theorem, without any relations on them. 

As the first step, we describe all $\gl(W)$-invariant maps $\hat{\phi}_0^{(N)}$, by the H.Weyl invariant theorem.

One easily sees that the vector space of such invariants has dimension $N+1$, and it can be described as follows. The first basic invariant $t_1$ is obtained as the post-composition with the symmetrization map
\begin{equation}\label{appa11}
W\to S^N(W)\xrightarrow{\mathrm{PBW}}T^N(W)
\end{equation}
and the remaining basic invariants $t_2,\dots,t_{N+1}$ are obtained as the composition(s):
\begin{equation}\label{appa12}
W^*\otimes S^N(W)\xrightarrow{\mathrm{div}}S^{N-1}(W)\xrightarrow{\mathrm{PBW}}T^{N-1}(W)\xrightarrow{e\in W^*\otimes W}
W^*\otimes T^N(W)
\end{equation}
where the rightmost map inserts the ``new'' factor $W$ after $\ell$ first factors, on the $(\ell+1)$-position, $\ell=0,\dots,N-1$, and $e\in W^*\otimes W$ is the canonical invariant.

Turn back to the statement of Lemma. If an imbedding $\phi_0$ existed, the components of the corresponding imbedding $\hat{\phi}_0^{(N)}$ were linear combinations of the described $N+1$ basic invariants,
\begin{equation}\label{appa13}
\hat{\phi}_0^{(N)}=a_{1N}t_1+a_{2N}t_2+\dots+a_{N+1,N}t_{N+1},\ \ a_{ij}\in\mathbb{C}
\end{equation}
and these components would define a map of Lie algebras, modulo the inner derivations of $T(W)$.

To get a contradiction, it is easier computing with divergence zero vector fields, so that $t_2=\dots=t_{N+1}=0$ on such vector fields.

That is, for {\it any} two divergence 0 (homogeneous) vector fields $v_1,v_2$ of degrees $M,N$ much smaller than $\dim W$, one would have
\begin{equation}\label{appa14}
[t_1(v_1),t_1(v_2)]=c_{MN}t_1(\{v_1,v_2\}) \ \text{\it modulo the inner derivations},\ c_{MN}\in\mathbb{C}
\end{equation}
where $t_1$ is the symmetrization map, see \eqref{appa11}.

We present two quadratic divergence 0 vector fields $v_1,v_2$ for which \eqref{appa14} fails.

Condider $\dim W\ge 3$, let $x,y,z$ be the first three coordinates. 

Take $$v_1=y^4\frac{\partial}{\partial x},\ v_2=
z^2\frac{\partial}{\partial y}$$

The commutator $$\{v_1,v_2\}= 4 z^2y^3 \frac{\partial}{\partial x}$$

On the other hand, $$[t_1(v_1),t(v_2)]=
(z^2 y^3+ y z^2 y^2+y^2 z^2 y +y^3 z^2)\frac{\partial}{\partial x}$$
and it is not equal to $$c\cdot t_1(4z^2y^3\frac{\partial}{\partial x})$$ even modulo the inner derivations. (The latter statement is true because all inner derivations in $\Der_\Assoc T(W)$ should contain $x$ as a factor in the coefficient at $\frac{\partial}{\partial x}$).

\end{proof}

Proposition \ref{propa2} is proven.
\end{proof}

Theorem \ref{theorema1} is proven.
\qed

\begin{remark}{\rm
If there existed a $\gl(W)$-equivariant isomorphism $\phi\colon \g_1\to \g_2$ of Lie algebras, the Kontsevich formality theorem would follow from it. It can be shown by taking the free Koszul resolution of the polynomial algebra $S(W)$. Let $R\to S(W)$ be this resolution. Then $R$ is a free associative dg algebra whose underlying graded vector space is $S(\Lie(S(W[1])[-1]))$. The differential in $R$ is tangential to $\Lie(S(W[1])[-1])$. Take $W_1=S(W[1])[-1]$. It is easy to see that $(\g_1(W_1), d_R)\simeq T_\poly(W^*)$, $(\g_2(W_1),d_R)\simeq \Hoch(S(W))$ as dg Lie algebras. This idea, suggested by Boris Feigin around '98-'99, was a source of inspiration for the author's work on this paper. 
}
\end{remark}

\begin{remark}{\rm
It would be very interesting to compute the Chevalley-Eilenberg cohomology $H_\CE^\udot(\g_1(W),\gl(W); \g_1(W))$ with the adjoint coefficients, at least in the inductive limit $\dim W\to\infty$. The argument from the previous Remark shows that the latter stable cohomology is closely related with the stable cohomology $H_\CE^\udot(T_\poly(W),\gl(W); T_\poly(W))$, which has been actively studied in the recent years, see e.g. [W].
}
\end{remark}
\section{\sc Where does Tamarkin's proof fail for an infinite-dimensional $V$}
Here we show where Tamarkin's proof of the Kontsevich formality for $S(V^*)$ does not work, where $\dim V=\infty$.
For the definitions of $S(V^*)$, $T_\fin(V)$, $\Hoch^\udot_\fin(S(V^*))$ see Sections \ref{section1.0}, \ref{section1.1}, and \ref{section1.2}.

The answer to the question in the title of this Appendix is the following. For a finite-dimensional vector space $W$, the graded commutative algebra $T_\poly(W^*)$ of polyvector fields is equal to $S(W\oplus W^*[-1])$ and therefore is {\it smooth} graded commutative algebra. Its smoothness is crucial in the proof of one of the key steps in the Tamarkin proof, see Proposition \ref{propb} below. Its smoothness implies the vanishing of higher Andre-Quillen cohomology.

On the other side, for an infinite-dimensional $V$, the graded commutative algebra $T_\fin(V^*)$ is not of the form $S(L)$, for a graded vector space $L$, as follows immediately from its definition. In fact, it fails to be smooth. It results in presence of non-trivial higher (Andre-Quillen-like) cohomology, which breaks the proof of Proposition \ref{propb}.

We provide more detail below.

The key-point of Tamarkin's proof of the Kontsvich formality is the following fact:
\begin{prop}\label{propb}
Let $W$ be a finite-dimensional vector space over a field of characteristic 0. Consider the polynomials polyvector fields $T_\poly(W)$ as a Gerstenhaber algebra (aka 2-algebra). Then the cohomology of the $\Aff(W)$-invariant part of the deformation complex $\mathrm{Def}_{\mathsf{hoe}_2}(T_\poly(W))$ of $T_\poly(W)$ as a homotopy 2-algebra, vanishes in all degrees. 
\end{prop}

We recall the main steps in the proof, because the failure of Tamarkin's proof for an infinite-dimensional vector space $V$ comes from the failure of this Proposition. 

Let $W$ be a finite-dimensional vector space. The deformation complex $\mathrm{Def}_{\mathsf{hoe}_2}(T_\poly(W))=\mathrm{Def}(\mathsf{e}_2\to\End_\mathrm{Op}(T_\poly(W)))$ is computed as the coderivations $\mathrm{Coder}(\Cobar_{\mathsf{e}_2}(T_\poly(W)))$. The Koszul dual cooperad to $\mathsf{e}_2$ is equal to $\mathsf{e}_2^*\{-2\}$, and $\Cobar_{\mathsf{e}_2}(X)=
(\mathsf{e}_2^*\{-2\}\circ X, d)$, where the differential $d$ comes from the $\mathsf{e}_2$-algebra structure on $X$. The differential $d=d_\Comm+d_{\Lie}$ has two components, expressed via the graded commutative product, and via the shifted by 1 Lie bracket, correspondingly. See ??? for more detail.

For any graded vector space $T$, one has
$$
\mathsf{e}_2\circ T=S^*(\Lie^*(T[-1])[1])
$$
(where $S^*(-)$ and $\Lie^*(-)$ stand for the cofree commutative (corresp., Lie) coalgebra),
and thus
$$
\mathsf{e}_2^*\{-2\}(T_\poly(W))=(\mathsf{e}_2^*\circ T_\poly(W)\{2\})\{-2\}=S^*(\Lie^*(T_\poly(W)[1])[1])[-2]
$$
The corresponding complex of coderivations is equal
$$
\begin{aligned}
\ &\mathrm{Def}_{\mathsf{hoe}_2}(T_\poly(W))=\Hom_\k\big(S^*(\Lie^*(T_\poly(W)[1])[1])[-2],\ T_\poly(W)\big)=\\
&
\Hom_\k\big(S^*(\Lie^*(T_\poly(W)[1])[1]),\ T_\poly(W)\big)[2]
\end{aligned}
$$
The differential has two components, $d=d_\Comm+d_\Lie$. One uses the spectral sequence of the bicomplex, and computes the cohomology of $d_\Comm$ at the first step.

One has:
\begin{equation}
\big(\mathrm{Def}_{\mathsf{hoe}_2}(T_\poly(W)),\ d_\Comm\big)=\prod_{k\ge 1}\Big(\Hom_\k(S^k(\Lie^*(T_\poly(W)[1])[1]),T_\poly)[2],\ d_\Comm\Big)
\end{equation}
(the component with $k=0$ does not figure out in the deformation complex).

Each component, corresponded to any $k\ge 1$, is a sub-complex with respect to the differential $d_\Comm$.
The simplest one among them is corresponded to the case $k=1$, it is 
\begin{equation}
C_\mathrm{Harr}(T_\poly(W))=\big(\Hom_\k(\Lie(T_\poly(W)[1],T_\poly(W))[1],d_\Comm\big)
\end{equation}
It is the Harrison complex, computing the Andre-Quillen cohomology.

It is well-known that {\it for any smooth over $\k$} (graded) commutative algebra $A$, the complex $C_\mathrm{Harr}(A)$ {\it has vanishing cohomology except for the degree 0}, where its cohomology is equal to $\Der_\Comm(A)$:
\begin{equation}\label{appb1}
H^\udot(C_\mathrm{Harr}(A))=\Der_\Comm(A)[0]
\end{equation}
Analogously, for higher $k$, one has
\begin{equation}\label{appb2}
H^\udot(\Hom_\k(S^k(\Lie(A[1])[1]),A[2]),d_\Comm))=S_{A}^k(\Der_\Comm(A)[-2])[2]
\end{equation}
The rest of the computation is not relevant for our goal in this Appendix, and we refer the interested reader to [T] for more detail on the Tamarkin computation.

\hspace{1mm}

The matter is that the {\it smoothness} of $A$ is essential for vanishing of the higher cohomology in \eqref{appb1} and \eqref{appb2}.

Turn back to our case of an invinite-dimensional vector space $V$, as in Section \ref{section1.0}.
Then the polyvector fields $T_\fin(V)$ {\it is not smooth as a commutative algebra}.
Indeed, to establish the smoothness of $T_\poly(W)$ for a finite-dimensional $W$, we use
\begin{equation}\label{appb3}
T_\poly(W)=S(W^*\oplus W[-1])
\end{equation}
which is smooth as a polynomial algebra.

For $T_\fin(V)$ (see Section \ref{section1.1} for the definition), one does not have an equality similar to \eqref{appb3}:
$$
T_\fin(V)\ne S(-)
$$
In fact, the graded commutative algebra $T_\fin(V)$ is not smooth. 

It implies that the higher  homology of the complexes $\Hom_\k((S^k(\Lie(A[1])[1]),A[2]),d_\Comm)$, $k\ge 1$, may not vanish. It results in the failure of the overall argument.

\bigskip

{\small
\noindent {\sc Universiteit Antwerpen, Campus Middelheim, Wiskunde en Informatica, Gebouw G\\
Middelheimlaan 1, 2020 Antwerpen, Belgi\"{e}}}

\bigskip

\noindent{{\it e-mail}: {\tt Boris.Shoikhet@uantwerpen.be}}

\end{document}